\documentclass[12pt, reqno]{amsart}
\usepackage{pifont}
\usepackage{mathrsfs}
\usepackage{geometry}
\usepackage{titletoc}
\usepackage{stix2}
\usepackage{amscd}

\usepackage{amsmath}
\usepackage{amssymb} 
\usepackage{enumitem} 
\usepackage{mathtools} 
\usepackage[table]{xcolor} 
\usepackage[all]{xy} 
\usepackage{tikz} 
\usepackage{tikz-cd}
\usepackage{indentfirst} 
\usepackage{babel} 
\usepackage{setspace} 

\usepackage[colorlinks,linkcolor=red,anchorcolor=green,citecolor=blue]{hyperref} 
\hypersetup{linktocpage = true} 

\usepackage{rotating} 

\usepackage{ytableau} 
\usepackage{longtable} 
\newcolumntype{M}[1]{>{\centering\arraybackslash}m{#1}} 

\usepackage{fancyhdr}

\newcommand\sE{{\mathscr E}}

\newcommand\sL{{\mathscr L}}

\newcommand\bp{{\bar\partial}}

\usepackage{amsthm}
\theoremstyle{plain}

\newtheorem{thm}{Theorem}[section]
\newtheorem{lemma}[thm]{Lemma}
\newtheorem{prop}[thm]{Proposition}
\newtheorem{cor}[thm]{Corollary}
\newtheorem{defn}[thm]{Definition}

\newtheorem{problem}[thm]{Problem}

\theoremstyle{definition}
\newtheorem{example}[thm]{Example}
\newtheorem{remark}[thm]{Remark}

\newcommand{\btheorem}{\begin{thm}}
    \newcommand{\etheorem}{\end{thm}}
\newcommand{\bproposition}{\begin{prop}}
    \newcommand{\eproposition}{\end{prop}}
\newcommand{\bdefinition}{\begin{defn}}
    \newcommand{\edefinition}{\end{defn}}
\newcommand{\bcorollary}{\begin{cor}}
    \newcommand{\ecorollary}{\end{cor}}
\newcommand{\bproof}{\begin{proof}}
    \newcommand{\eproof}{\end{proof}}
\newcommand{\bremark}{\begin{remark}}
    \newcommand{\eremark}{\end{remark}}
\newcommand{\eexample}{\end{example}}
\newcommand{\bexample}{\begin{example}}

\newcommand{\elemma}{\end{lemma}}
\newcommand{\blemma}{\begin{lemma}}

\newcommand{\sq}{\sqrt{-1}}

\newcommand{\p}{\partial}

\renewcommand{\bar}{\overline}

\renewcommand{\phi}{\varphi}

\newcommand{\beq}{\begin{equation}}
\newcommand{\eeq}{\end{equation}}
\newcommand{\ee}{\end{eqnarray*}}
\newcommand{\be}{\begin{eqnarray*}}

\newcommand{\bd}{\begin{enumerate}}
    \newcommand{\ed}{\end{enumerate}}

\renewcommand{\tilde}{\widetilde}

\newcommand{\qtq}[1]{\quad\mbox{#1}\quad}
\renewcommand{\bp}{\bar{\partial}}
\newcommand{\Om}{\Omega}

\newcommand{\ts}{\otimes}

\renewcommand{\>}{\rightarrow}

\newcommand{\C}{{\mathbb C}}

\newcommand{\R}{{\mathbb R}}

\renewcommand{\>}{\rightarrow}

\renewcommand{\p}{{\partial}}
\renewcommand{\bp}{{\bar{\partial}}}

\newcommand{\vone}{ \vskip 1\baselineskip}
\newcommand{\om}{\omega}

\renewcommand{\bar}{\overline}
\renewcommand{\tilde}{\widetilde}

\newcommand{\smo}{\sqrt{-1}}

\newcommand{\tr}{\mathrm{tr}}
\newcommand{\nm}[1]{\left\Vert #1\right\Vert}
\setlist[itemize]{leftmargin=*}
\setlist[enumerate]{leftmargin=*}

\numberwithin{equation}{section} 

\makeatletter

\usepackage{fancyhdr}
\title{Existence of Hermitian metrics  with  prescribed Hermitian-Yang-Mills tensors  II}

\author{Jiaxuan Fan}
\author{Mingwei Wang}
\author{Xiaokui Yang}
\author{Shing-Tung Yau}

\address{Jiaxuan Fan, Qiuzhen College, Tsinghua University, Beijing, 100084, China}
\email{fanjiaxu21@mails.tsinghua.edu.cn}

\address{Mingwei Wang, Qiuzhen College, Tsinghua University, Beijing, 100084, China}
\email{wangmw21@mails.tsinghua.edu.cn}

\address{Xiaokui Yang, Department of Mathematics and Yau Mathematical Sciences Center, Tsinghua University, Beijing, 100084, China}
\email{xkyang@mail.tsinghua.edu.cn}

\address{Shing-Tung Yau, Yau Mathematical Sciences Center and  Qiuzhen College, Tsinghua University, Beijing, 100084, China}
\email{styau@mail.tsinghua.edu.cn}

\begin{document}

    \begin{abstract} In this paper, we solve the prescribed Hermitian-Yang-Mills tensor problem for Higgs bundles over compact complex manifolds.
        Let $ (E,\theta) $ be a Higgs  bundle over a compact Hermitian manifold $(M,\omega_g) $. Suppose that there exists a  Hermitian metric $ h_0 $ on $E$ such that the Hermitian-Yang-Mills tensor $ \Lambda_{\omega_g}\left(\sq  R^{D^{h_0}}\right) $ of the Higgs connection   is  positive-definite. Then for any positive-definite Hermitian  tensor $ P\in \Gamma\left(M,E^*\otimes \bar E^*\right) $, there exists  a unique  Hermitian metric $ h $ on $E$ such that  $$\Lambda_{\omega_g} \left(\sq R^{D^h}\right)=P.$$
    We also establish quantitative Chern number inequalities for Higgs bundles.

    \end{abstract}

    \maketitle {
        \setcounter{tocdepth}{1}

    {\small{    \begin{spacing}{1.1} \tableofcontents %
                \dottedcontents{section}[1.8cm]{}{3em}{5pt} %
\end{spacing} }} }

    \section{Introduction}

    This paper continues the investigation initiated in \cite{WYY26+}, where we solved the prescribed Hermitian-Yang-Mills tensor equation for holomorphic vector bundles on compact K\"ahler manifolds, which is a vector bundle version of the Calabi-Yau theorem.
    Let's recall the set-up briefly.  On a compact K\"ahler manifold $(M,\omega_g)$, it is well-known that the class of Ricci curvature $\mathrm{Ric}(\omega_g)$ represents the first Chern class $c_1(M)$ up to a factor of $2\pi$, i.e.,
    $$\left[\mathrm{Ric}(\omega_g)\right]=2\pi c_1(M)\in H^{1,1}(M,\R).$$
    The celebrated solution to the Calabi conjecture, proven by S.-T. Yau \cite{Cal57,Yau78}, establishes the converse: for any closed real $(1,1)$ form $\Om$ representing $2\pi c_1(M)$,  there exists a unique K\"ahler metric $\omega\in [\omega_g]$ such that
    \beq \mathrm{Ric}(\omega)=\Om. \label{PRicci}\eeq
    The Calabi-Yau theorem has indeed been extended beyond K\"ahler manifolds to more general Hermitian manifolds.  We refer to \cite{FY08, TW10a, TW10b, STW17} and the references therein.\\

        The classical Donaldson-Uhlenbeck-Yau theorem (\cite{Don85}, \cite{UY86}, \cite{Don87}, see also \cite{NS65}),  establishes a deep equivalence on compact K\"ahler manifolds: a holomorphic vector bundle admits a Hermitian-Einstein metric if and only if it is polystable in the sense of Mumford-Takemoto. This result forged a critical link between differential geometry (the existence of special metrics solving a non-linear PDE) and algebraic geometry (a purely algebraic condition of stability). C. Simpson's pivotal insight was to generalize this correspondence to the setting of Higgs bundles and demonstrated that the Donaldson-Uhlenbeck-Yau correspondence could be extended to this context.
     It is well-known that the concept of Higgs bundles was first introduced by N. Hitchin in \cite{Hit87}, in the context of his  work on self-duality equations over Riemann surfaces.  C. Simpson subsequently generalized Hitchin's  construction to higher-dimensional K\"ahler manifolds, and proved that a Higgs bundle over a compact K\"ahler manifold admits a Hermitian-Einstein metric compatible with the Higgs structure if and only if it is Higgs polystable. This work lays the foundational groundwork for the Hitchin-Simpson correspondence--a critical bridge linking Higgs bundles and local systems. It illuminates deep interdisciplinary connections spanning algebraic geometry, differential geometry, and representation theory.  For comprehensive technical details and contextual background, we refer the reader to Simpson's seminal works \cite{Sim88, Sim92, Sim94a, Sim94b} and the extensive references cited therein.\\

        Let $M$ be a compact complex manifold and $E$ be a holomorphic vector bundle over $M$.   Let $\theta\in \Om^{1,0}(M,E^*\ts E)$ be an $ E^* \otimes E $-valued smooth $(1,0)$-form.  $(E, \theta)$ is called a Higgs bundle if $\bp\theta=0$.  Suppose that $h$ is a smooth Hermitian metric on $E$. The Hermitian adjoint of $\theta$ with respect to $h$ is denoted by $\theta^\star_h$. The Higgs connection $D^{h}$ for the Hermitian Higgs bundle $(E,h,\theta)$ is defined as
        \beq  D^h=\nabla^h+\theta+\theta_h^\star\eeq
    where $\nabla^h$ is the Chern connection of $(E,h)$. The curvature tensor  of the Higgs connection $D^h$ is denoted by \beq R^{D^h}:=(D^h)^2\in \Gamma(M,\Lambda^2 T^*M\ts E^*\ts E).\eeq
A Higgs bundle is called Hermitian-Einstein if there exist a constant $\lambda_0\in \R$ and a smooth Hermitian metric $h$ on $E$ such that the equation
\beq \Lambda_{\omega_g}\left(\sq R^{D^h}\right)=\lambda_0\cdot  h \label{HYMH}\eeq
holds as an equality of tensors in $\Gamma(M, E^*\ts \bar E^*)$. It is proved in \cite{Sim88} that \eqref{HYMH} is solvable if and only if $(E,\theta)$ is Higgs polystable over a compact K\"ahler manifold. Note that, if $\theta=0$, the equation \eqref{HYMH} is the classical Hermitian-Einstein equation.
 In what follows, we refer to the tensor $\Lambda_{\omega_g}\left(\sq R^{D^h}\right)$ as the \emph{Hermitian-Yang-Mills-Higgs tensor}.  \\

The main result of this paper is the following analogue of the Calabi-Yau theorem, addressing the prescribed Hermitian-Yang-Mills-Higgs tensor problem for Higgs bundles on compact complex manifolds.
        \btheorem
    \label{main1}
    Let $(E,\theta)$ be a Higgs bundle over a compact Hermitian manifold $(M,\omega_g) $. Suppose that there exists a  Hermitian metric $ h_0 $ on $E$ such that the Hermitian-Yang-Mills-Higgs tensor of  $(E,h_0,\theta)$ satisfies $$ \Lambda_{\omega_g}\left(\sq  R^{D^{h_0}}\right)>0.$$ Then for any positive-definite Hermitian  tensor $ P\in \Gamma(M,E^*\otimes \bar E^*) $, there exists a unique Hermitian metric $ h $ on $E$ such that  \beq \Lambda_{\omega_g}\left(\sq R^{D^h}\right)=P.\eeq
    \etheorem
    \noindent The positivity of  $ \Lambda_{\omega_g}\left(\sq  R^{D^{h_0}}\right)$ can be weakened to the positivity of the integral of the minimum eigenvalue of
    $\Lambda_{\omega_g}\left(\sq  R^{D^{h_0}}\right) $ over a Gauduchon manifold (Theorem \ref{main6}).
    \noindent  The following special case--where $\theta=0$--merits particular attention:
    \bcorollary
    \label{main2}
    Let $E$ be a holomorphic vector bundle over a compact Hermitian manifold $(M,\omega_g)$. Suppose that there exists a Hermitian metric $ h_0 $ on $E$ such that the Hermitian-Yang-Mills tensor of the Chern connection satisfies $$ \Lambda_{\omega_g}\left(\sq  R^{{h_0}}\right)>0.$$ Then for any  positive-definite Hermitian tensor $ P\in \Gamma(M,E^*\otimes \bar E^*) $, there exists a unique Hermitian metric $ h $ on $E$ such that  \beq \Lambda_{\omega_g}\left(\sq R^{h}\right)=P.\eeq
    \ecorollary

        \noindent The uniqueness part of Theorem \ref{main1} is established by the following comparison principle for Hermitian-Yang-Mills-Higgs tensors:

    \btheorem   \label{main3}
    Let $(E,\theta)$ be a Higgs bundle over a compact Hermitian manifold $(M,\omega_g)$. Suppose that $h$ and $h_0$ are Hermitian metrics on $E$. If the Hermitian-Yang-Mills tensors of their Higgs connections satisfy $
    \Lambda_{\omega_g}\left(\smo R^{D^{h_0}}\right) > 0 $ and \beq \Lambda_{\omega_g}\left(\smo R^{D^h}\right) \leq \Lambda_{\omega_g}\left(\smo R^{D^{h_0}}\right)  \label{Hermitiantensor}\eeq as Hermitian tensors in $ \Gamma(M,E^*\otimes \bar E{}^*) $, then $ h \leq h_0 $.
    \etheorem

\noindent   Note that inequality \eqref{Hermitiantensor} holds as Hermitian tensors in $\Gamma(M,E^*\otimes  \bar E^*)$,  \textbf{but not  as scaling invariant endomorphisms} in  $\Gamma(M,E^*\otimes E)$.
 \noindent  When $\theta=0$, we obtain:
    \bcorollary\label{main4} Let $E$ be a holomorphic vector bundle over a compact Hermitian manifold $(M,\omega_g)$. Suppose that $h$ and $h_0$ are Hermitian metrics on $E$. If the Hermitian-Yang-Mills tensors of their Chern connections satisfy $
    \Lambda_{\omega_g}\left(\smo R^{h_0}\right) > 0 $ and \beq \Lambda_{\omega_g}\left(\smo R^{h}\right) \leq \Lambda_{\omega_g}\left(\smo R^{{h_0}}\right)  \eeq as Hermitian tensors in $ \Gamma(M,E^*\otimes \bar E{}^*) $, then $ h \leq h_0 $.
    \ecorollary

        \noindent The K\"ahler case of Corollary \ref{main2} was resolved in \cite{WYY26+}.  We adopt  approaches similar to those in \cite{WYY26+}. Let us clarify the core challenges and distinguishing aspects in proving Theorem \ref{main1} and Theorem \ref{main3}. First, the Higgs connection typically fails to be metric-compatible, which complicates the control of curvature terms. Second, both the torsion of the Hermitian manifold and the Higgs field  $\theta$ introduce  analytical challenges, requiring refined uniform estimates to control their global behavior.  Last but not least, the linearized equations deviate substantially from those in the K\"ahler case, primarily due to the lack of K\"ahler identities, requiring a more delicate treatment on the injectivity of the associated elliptic systems (e.g. Theorem \ref{injectivity}).\\

    Motivated by the comparison result established in Theorem \ref{main3}, we demonstrate a quantitative form of Chern number inequalities, extending those known for Hermitian-Einstein metrics on stable Higgs bundles (\cite{Sim88}).

    \btheorem \label{main5} Let $(M,\omega_g)$ be a compact K\"ahler manifold and $(E,\theta)$ be an integrable Higgs vector bundle of rank $r$ over $M$.  Suppose that there exists a smooth Hermitian metric $h_0$ on $E$ such that the Hermitian-Yang-Mills tensor of the Higgs connection satisfies
    \beq  a\cdot h_0 \leq   \Lambda_{\omega_g}\left(\smo  R^{D^{h_0}}\right) \leq b \cdot h_0, \eeq
    for some constants $a, b\in \R$. Then the following Chern number inequality holds
    \beq \int_M \left( (r-1)c^2_1(E)-2rc_2(E) \right) \wedge \omega_g^{n-2} \leq \frac{r(r-1)\left(b - a\right)^2}{8\pi^2n^2}\int_M \omega_g^n.\label{CN} \eeq
    \etheorem
        \noindent\textbf{Acknowledgements}. The third named author would like to thank  Bing-Long Chen, Jixiang Fu,  Kefeng Liu and   Valentino Tosatti  for inspiring  discussions.

    \vskip 1\baselineskip

    \section{Background materials}
In this section, we provide a concise overview of foundational concepts related to Higgs bundles. For more exhaustive, in-depth treatments of the subject, we direct readers to the seminal works \cite{Hit87}, \cite{Sim88}, \cite{Sim92} and \cite{Sim94a}. Let $M$ be a compact complex manifold and $E$ be a holomorphic vector bundle over $M$.   Let $\theta\in \Om^{1,0}(M,E^*\ts E)$ be an $ E^* \otimes E $-valued smooth $(1,0)$-form.  Recall that
\bd \item  If $\bp \theta=0$,  then $(E,\theta) $ is called a \emph{Higgs bundle}.
\item If $\bp\theta=0$ and $\theta\wedge\theta=0$,  $(E,\theta)$ is called an \emph{integrable Higgs bundle}.
\ed
There is a natural operator $ D^{\prime\prime}: \Gamma(M, E) \> \Omega^1(M,E)$
given by
\beq D^{\prime\prime}v = \bar\p v + \theta v. \eeq
For any smooth Hermitian metric $ h_0 $ on $ E $, let $\nabla^{h_0}$ be the Chern connection of $(E,h_0)$. There is a decomposition $\nabla^{h_0}=\p^{h_0}+\bp$ where $ \p^{h_0} $ is the $(1,0)$-part of the Chern connection.  The Hermitian adjoint of the Higgs field $\theta$ with respect to $h_0$ is denoted by $\theta_{h_0}^{\star}$, i.e., for any $v, w\in \Gamma(M,E)$,
\beq h_0(\theta v,w )= h_0\left(v, \theta_{h_0}^{\star}w\right). \eeq
In local holomorphic coordinates $\{z^i\}$ of $M$ and local holomorphic basis $ \{ e_\alpha \} $ of $E$, if we write $\theta=\theta_{i\alpha}^\beta dz^i\ts e^\alpha\ts e_\beta$, then $\theta_{h_0}^{\star}=(\theta_{h_0}^{\star})_{\bar j \gamma}^\delta d\bar z^j \ts e^\gamma\ts e_\delta\in\Om^{0,1}(M,E^*\ts E)$ where
\beq \left(\theta_{h_0}^{\star}\right)_{\bar j\gamma}^\delta = h_{0,\gamma\bar\beta}h_0^{\delta\bar\mu} \; \overline{\theta_{j\mu}^\beta}.\label{conjugate}
\eeq
Let  $ D^{\prime,h_0}:\Gamma(M,E) \> \Omega^1(M,E) $  be the operator defined by
\beq D^{\prime,h_0} v = \p^{h_0}v + \theta_{h_0}^\star v. \eeq
We call the affine connection  \beq D^{h_0}=D^{\prime,h_0}+D'' \eeq   the Higgs connection on $(E,h_0,\theta)$.
In general, it is not a metric compatible connection.  The curvature tensor $\Theta^{D^{h_0}}\in \Gamma(M,\Lambda^2T^*M\ts E^*\otimes E)$ of this connection is defined by
\beq \Theta^{D^{h_0}} :=(D^{h_0})^2=(D'')^2+D^{\prime,h_0}D^{\prime\prime} + D^{\prime\prime}D^{\prime,h_0}+(D'^{,h_0})^2 . \eeq
By using the Hermitian metric $h_0$, we can define
\beq R^{D^{h_0}}:=\Theta^{D^{h_0}}\cdot h_0 \in \Gamma(M,\Lambda^2T^*M\ts E^*\otimes \bar E^*)\eeq
The \textbf{Hermitian-Yang-Mills tensor} $ S^{h_0}\in
\Gamma(M,E^*\ts \bar E^*)$ of $(E,h_0,\theta)$  is defined as \beq S^{h_0}
:=\Lambda_{\omega_g}\left(\smo R^{D^{h_0}}\right)=\left(g^{i\bar
    j}R^{D^{h_0}}_{i\bar j\alpha\bar\beta}\right) e^\alpha\ts \bar e^\beta \in
\Gamma(M,E^*\ts \bar E^*).
\eeq
We say that $S^{h_0}$ is positive-definite, $S^{h_0}>0$, if $\left(S^{h_0}_{\alpha\bar\beta}\right)=\left(g^{i\bar
    j}R^{D^{h_0}}_{i\bar j\alpha\bar\beta}\right)$ is a positive-definite Hermitian matrix at each point of $M$. We also set \beq   K^{h_0}:=S^{h_0}\cdot h_0^{-1}=\left(g^{i\bar j}\left(R^{D^{h_0}}\right)_{i\bar j\alpha}^{\beta}\right) e^\alpha\ts  e_\beta\in \Gamma(M,E^*\ts E).
\eeq
The formal adjoints of $D^{\prime, h_0}$ and $D^{\prime\prime}$ with respect to Hermitian metrics $g$ and $h_0$ are denoted by $(D^{\prime, h_0})^*$ and $(D^{\prime\prime})^*$ respectively.  These operators can be extended to $\Om^{p}(M,E)$ in the  standard manner.\\

Recall that for $ P=P_\alpha^\beta e^\alpha\ts e_\beta \in \Omega^p(M,E^*\otimes E) $ and $ Q=Q_\alpha^\beta e^\alpha\ts e_\beta \in \Omega^q(M,E^*\otimes E) $, the product \beq  P\cdot Q \in \Omega^{p+q}(M,E^*\otimes E) \eeq  is defined as
\beq P \cdot Q = (P_\alpha^\beta \wedge Q_\beta^\gamma) \; e^\alpha \otimes e_\gamma. \eeq

Let $ h $ be  another smooth Hermitian metric on $ E $.  We set $ H = h \cdot h_0^{-1} \in \Gamma(M,E^*\otimes E)$.  It is easy to see that
\beq
H=h_{\alpha\bar{\gamma}}h_0^{\beta\gamma}e^\alpha\otimes e_{\beta}
\eeq
Moreover, for any $s,t\in \Gamma(M,E)$,
\beq h_0(Hs,t)=h(s,t).\label{H}
\eeq

\noindent There following linear algebraic calculation rules are given in \cite[Lemma~2.2]{WYY26+}.
\blemma \label{linearalgebra}\bd
\item $H=h\cdot h_0^{-1}\in\Gamma(M,E^*\ts E)$    is  $h_0$-Hermitian.

\item If $P\in \Gamma(M,E^*\ts E)$ is $h_0$-Hermitian,  for any $ A, B \in \Gamma(M,E^*\otimes E) $,
\beq h_0(P\cdot A, B)=h_0(A, P\cdot B), \quad  h_0(A \cdot P, B) = h_0(A, B\cdot P). \eeq

\item If $ A, B \in \Gamma(M,E^*\otimes E) $ are $h_0$-Hermitian and $ A\geq 0, B \geq 0 $, then
\beq h_0(A,B) \geq 0. \eeq Moreover,  for  any section $ C \in \Gamma(M,E^* \otimes E) $, one has
\beq h_0(A \cdot C \cdot B, C ) \geq 0. \eeq

\ed
\elemma

\noindent The Higgs connection $D^{h_0}$ can be extended to $\Om^{p}(M,E^*\ts E)$: for any $P\in \Om^{p}(M,E^*\ts E)$
\beq  D'' P=\bp P+\theta(P), \quad  D'^{,h_0}(P)=\p^{h_0}P+\theta^\star_{h_0}(P). \eeq

\noindent
The following computational formulas will be applied frequently throughout this work.
\blemma\label{keycomputation}
For any $P\in \Om^p(M,E^*\otimes E)$,
\beq
D''P=\bp P+(-1)^p P\cdot\theta-\theta\cdot P,\quad
D'^{,h_0}P=\p^{h_0}P+(-1)^p P\cdot \theta_{h_0}^\star-\theta_{h_0}^\star\cdot P. \label{connection1}
\eeq
In particular,
\beq \theta(P)=(-1)^p P\cdot\theta-\theta\cdot P, \quad \theta^\star_{h_0}(P)=(-1)^p P\cdot \theta_{h_0}^\star-\theta_{h_0}^\star\cdot P.\eeq
\elemma
\bproof Suppose that $P=P_\alpha^\beta e^\alpha\ts e_\beta$ where $P_{\alpha}^\beta$ are $p$-forms. It is obvious that
\be  D^{\prime\prime} P&=&\bp(P_{\alpha}^\beta) e^\alpha\ts e_\beta+(-1)^p P_{\alpha}^\beta\wedge  \left(D'' e^\alpha\right)\ts e_\beta + (-1)^p P_{\alpha}^\beta\wedge   e^\alpha\ts \left(D''e_\beta\right)\\
&=& \bp P-(-1)^p P_{\alpha}^\beta\wedge  \theta_{i\gamma}^\alpha dz^i\ts e^\gamma\ts e_\beta+(-1)^p P_{\alpha}^\beta\wedge  \theta_{i\beta}^\gamma dz^i \ts e^\alpha\ts e_\gamma.
\ee
It is clear that $(-1)^p P_{\alpha}^\beta\wedge  \theta_{i\gamma}^\alpha dz^i\ts e^\gamma\ts e_\beta=\theta\cdot P$ and $P_{\alpha}^\beta\wedge  \theta_{i\beta}^\gamma dz^i \ts e^\alpha\ts e_\gamma=P\cdot \theta$. The second formula admits a similar proof.
\eproof

\blemma
The following identities hold:
\beq \p^{h_0}\theta_{h_0}^\star=0, \quad
\theta_{h}^{\star}=H\cdot \theta_{h_0}^{\star}\cdot H^{-1}.\label{relation2}\eeq
Moreover, for any $s\in \Gamma(M,E)$, \beq
D'^{,h}(s)=H^{-1}\left( D'^{,h_0}(Hs)\right).\label{relation1}
\eeq
\elemma
\bproof  Let $A\in \Om^1(M,E^*\ts E)$ and $A^\star$ be the Hermitian adjoint of $A$ with respect to $h_0$.
Since the Chern connection is metric compatible, for any $s,t\in \Gamma(M,E)$,
\begin{align*}
h_0\left(\left(\p^{h_0}A\right) s,t\right)=&h_0\left(\p^{h_0}\left(A s\right),t\right)+h_0\left(A(\p^{h_0} s),t\right)\\
=&\p h_0(A s,t)+h_0(A s,\bp t)-h_0(\p^{h_0}s,A^\star t)\\
=&\p h_0(s,A^\star t)+h_0\left(s,A^\star(\bp t)\right)-h_0(\p^{h_0}s,A^\star t)\\
=& h_0\left(s,\bp(A^\star t)\right)+h_0\left(s,A^\star(\bp t)\right)=h_0\left(s,\left(\bp A^\star\right) t\right).
\end{align*}
Hence, we obtain the relation
\beq (\p^{h_0} A)^\star=\bp A^\star.\label{adjoint}\eeq
In particular, $\p^{h_0}\theta_{h_0}^\star=(\bp\theta)^\star=0$. Moreover,
by definition \eqref{conjugate},
\begin{align*}
H\cdot \theta_{h_0}^{\star}\cdot H^{-1}=&\bar{\theta_{i\delta}^\gamma}h_{0,\mu\bar{\gamma}}h_0^{\nu\bar{\delta}}H^\mu_\alpha (H^{-1})_\nu^\beta d\bar{z}^i\otimes e^\alpha\otimes e_{\beta}\\
=&\bar{\theta_{i\delta}^\gamma}h_0^{\mu\bar{\eta}}h_{0,\mu\bar{\gamma}}h_0^{\nu\bar{\delta}}h_{0,\nu\bar{\tau}} h_{\alpha\bar{\eta}}  h^{\beta\bar{\tau}} d\bar{z}^i\otimes e^\alpha\otimes e_{\beta}\\
=&\bar{\theta_{i\delta}^\gamma}h_{\alpha\bar{\gamma}}  h^{\beta\bar{\delta}} d\bar{z}^i\otimes e^\alpha\otimes e_{\beta}=\theta_{h}^{\star}.
\end{align*}
Hence, we obtain \eqref{relation2}.
On the other hand, for any $s,t\in \Gamma(M,E)$,
\begin{align}
\p h_0(Hs,t)=&h_0(\p^{h_0}(Hs),t)+h_0(Hs,\bp t),\\
\p h(s,t)=&h(\p^{h}s,t)+h(s,\bp t).
\end{align}
By using \eqref{H}, one  obtains
\beq
\p^{h_0}(Hs)=H(\p^{h}s).
\eeq
That is $\p^{h_0}(Hs)=H(\p^{h}s)$. By \eqref{relation2}, one has
\beq \theta_{h}^{\star}(s)=H^{-1}\left( \theta_{h_0}^{\star}\left(Hs\right) \right).\eeq
The summation yields the identity given in \eqref{relation1}.
\eproof

\blemma\label{Lebniz}
For any $ A, B \in \Gamma(M,E^*\otimes E) $, the following identities  hold:
\beq D'^{,h_0}(A\cdot B) = D'^{,h_0}A \cdot B + A \cdot D'^{,h_0}B, \quad  D''(A\cdot B) = D''A \cdot B + A \cdot D''B. \eeq
\elemma

\bproof It is easy to see that
\beq \p^{h_0}(A\cdot B) = \p^{h_0}A \cdot B + A \cdot \p^{h_0}B.\label{productderivative} \eeq
By \eqref{connection1},
\begin{align*}
D'^{,h_0}A \cdot B + A \cdot D'^{,h_0}B=&\p^{h_0}A\cdot B+A\cdot \theta^\star_{h_0}\cdot B-\theta^\star_{h_0}\cdot A\cdot B\\
&+A\cdot \p^{h_0}B+A\cdot B\cdot \theta^\star_{h_0}-A\cdot\theta^\star_{h_0}\cdot B\\
=&\p^{h_0}(A\cdot B)+A\cdot B\cdot \theta^\star_{h_0}-\theta^\star_{h_0}\cdot A\cdot B\\
=&D'^{,h_0}(A\cdot B).
\end{align*}
The other identity can be established in a similar way.
\eproof

\noindent The following refined Bochner-Kodaira formulas hold on Higgs bundles:
\blemma Let $(M,\omega_g)$ be a compact Hermitian manifold and $(E,h_0,\theta)$ be a Higgs bundle.
The following identities hold on $\Om^{p,q}(M,E)$:
\beq
(D'^{,h_0})^*+\tau^*=\sqrt{-1}[\Lambda_{\om_g},D''], \quad
(D'')^*+\bar{\tau}^*=-\sqrt{-1}[\Lambda_{\om_g},D'^{,h_0}], \label{bkformula}
\eeq
where $\tau=[\Lambda,\p\om]$ is the torsion tensor of $\omega_g$ and $\tau^*$ is the formal adjoint of $\tau$.
\elemma

\bproof The following Bochner-Kodaira formula is well-known  on a compact Hermitian manifold $(M,\om_g)$ (e.g. \cite[Lemma~4.3]{LY12}):
\beq
(\p^{h_0})^*+\tau^*=\sqrt{-1}[\Lambda_{\om_g},\bp].
\eeq
Hence,  we only need to show the following identity
\beq
(\theta_{h_0}^\star)^*=\sqrt{-1}[\Lambda_{\om_g},\theta] \label{BK1}
\eeq
Let $I_{i}$ and $I_{\bar j}$ be the contraction operators $i_{\frac{\p}{\p z^i}}$ and $i_{\frac{\p}{\p \bar z^j}}$.  It is easy to compute that
\begin{align*}
\sqrt{-1}[\Lambda_{\om_g},\theta]=&g^{i\bar{j}}I_{\bar{j}}I_idz^k\wedge\theta_{k}-g^{i\bar{j}}dz^k\wedge I_{\bar{j}}I_i\theta_{k}\\
=&g^{i\bar{j}}I_{\bar{j}} \left( I_{i}dz^k\wedge+dz^k\wedge I_{i}\right)  \theta_k=g^{i\bar{j}}I_{\bar{j}}\theta_{i}.
\end{align*}
Here we use the fact $I_{i}dz^k\wedge+dz^k\wedge I_{i}=\delta_{i}^{k}$ in the last identity. On the other hand, for any $s\in \Om^p(M,E)$ and $t\in \Om^{p+1}(M,E) $, one has
\beq
(\theta_{h_0}^\star s,t)=\left(d\bar{z}^i\wedge \left(\theta_{h_0}^\star\right)_{\bar{i}}s,t\right)
=\left(\left(\theta_{h_0}^\star\right)_{\bar{i}}s,g^{i\bar{j}}I_{\bar{j}}t\right)
=(s,g^{i\bar{j}}I_{\bar{j}}\theta_{i}t).
\eeq Hence, \eqref{BK1} holds. The other identity can be shown in a similar way.
\eproof

\noindent
The following lemma explicates how the Hermitian-Yang-Mills tensors of the Chern connection differs from that of the Higgs connection:
\blemma

The following identity holds as Hermitian tensors  in $ \Gamma(M,E^*\otimes E) $:
\beq \Lambda_{\omega_g} \left(\smo \Theta^{D^{h_0}}\right) - \Lambda_{\omega_g}\left(\smo \Theta^{h_0}\right) =- \Lambda_{\omega_g}\smo(\theta \cdot  \theta_{h_0}^{\star}+ \theta_{h_0}^{\star}\cdot \theta).\label{chernhiggs} \eeq
\elemma
\bproof
Since $
D''=\bp+\theta$ and $  D'^{,h_0}=\p^{h_0}+\theta_{h_0}^{\star}$,
one has
$$
\Lambda_{\omega_g}\left(\smo \Theta^{D^{h_0}}\right)=\Lambda_{\omega_g}\smo\left(\bp\p^{h_0}+\p^{h_0}\bp+\bp\circ \theta+\theta\circ \bp +\p^{h_0}\circ \theta_{h_0}^{\star}+\theta_{h_0}^\star\circ \p^{h_0}+\theta\circ \theta_{h_0}^{\star}+\theta_{h_0}^{\star}\circ \theta \right).
$$
On the other hand, it is easy to see that
$$ \Lambda_{\omega_g}\smo\left(\bp\p^{h_0}+\p^{h_0}\bp\right)=\Lambda_{\omega_g}\left(\smo \Theta^{h_0}\right),$$
where $\Theta^{h_0}$ is the curvature tensor of the Chern connection on $(E,h_0)$.
Moreover, for any $s\in \Gamma(M,E)$,
\beq (\bp\circ \theta+\theta\circ \bp )s= (\bp\theta)(s)=0,\quad \left(\p^{h_0}\circ \theta_{h_0}^{\star}+\theta_{h_0}^\star\circ \p^{h_0}\right)(s)=(\p^{h_0}\theta_{h_0}^\star)(s)=0,\eeq
and
\beq \left(\theta\circ \theta_{h_0}^{\star}+\theta_{h_0}^{\star}\circ \theta \right)(s)=-(\theta \cdot  \theta_{h_0}^{\star}+ \theta_{h_0}^{\star}\cdot \theta)(s).\eeq
This completes the proof.
\eproof

\bcorollary\label{hermitian1}
$\Lambda_{\omega_g}\left(\sq \Theta^{D^{h_0}}\right)\in \Gamma(M,E^*\otimes E)$ is $h_0$-Hermitian.
\ecorollary

\noindent The following lemma describes the difference between the Hermitian-Yang-Mills tensors of two different Higgs connections.
\bproposition The following identity holds on $ \Gamma(M,E^*\otimes E) $:
\beq \Lambda_{\omega_g} \left(\smo \Theta^{D^h}\right) -\Lambda_{\omega_g} \left(\smo \Theta^{D^{h_0}}\right) = \Lambda_{\omega_g}\smo D^{\prime\prime}\left(D^{\prime,h_0}H \cdot H^{-1}\right).\label{curvaturedifference} \eeq
\eproposition
\bproof  We establish in the proof of \cite[Proposition~3.6]{WYY26+} the following formula on Hermitian manifolds for curvature of Chern connections:
\beq \Lambda_{\omega_g} \left(\smo \Theta^{h}\right) -\Lambda_{\omega_g} \left(\smo \Theta^{h_0}\right)= \Lambda_{\omega_g}\smo\bp\left( (\p^{h_0} H) \cdot H^{-1} \right).\label{difference1}
\eeq
By \eqref{chernhiggs} and \eqref{relation2}, one has
\beq
\Lambda_{\omega_g} \left(\smo \Theta^{D^{h_0}}\right) - \Lambda_{\omega_g}\left(\smo \Theta^{h_0}\right)=- \Lambda_{\omega_g}\smo(\theta \cdot  \theta_{h_0}^{\star}+ \theta_{h_0}^{\star}\cdot \theta)\label{difference2}
\eeq
and
\begin{eqnarray}
\Lambda_{\omega_g} \left(\smo \Theta^{D^{h}}\right) - \Lambda_{\omega_g}\left(\smo \Theta^{h}\right)&=&- \Lambda_{\omega_g}\smo(\theta \cdot  \theta_{h}^{\star}+ \theta_{h}^{\star}\cdot \theta) \label{difference3}\\
&=&- \Lambda_{\omega_g}\smo(\theta \cdot  H\cdot \theta_{h_0}^{\star}\cdot H^{-1}+ H\cdot \theta_{h_0}^{\star}\cdot H^{-1}\cdot \theta).\nonumber
\end{eqnarray}
By \eqref{difference1}, \eqref{difference2} and \eqref{difference3}, one obtains
\begin{align}
\begin{split}
&\Lambda_{\omega_g} \left(\smo \Theta^{D^h}\right) -\Lambda_{\omega_g} \left(\smo \Theta^{D^{h_0}}\right)-\Lambda_{\omega_g}\smo\bp\left( (\p^{h_0} H) \cdot H^{-1} \right)\\
=&\Lambda_{\omega_g}\smo\left(\theta \cdot  \theta_{h_0}^{\star}+ \theta_{h_0}^{\star}\cdot \theta-\theta \cdot  H\cdot \theta_{h_0}^{\star}\cdot H^{-1}- H\cdot \theta_{h_0}^{\star}\cdot H^{-1}\cdot \theta \right).\\
\end{split}
\label{difference4}
\end{align}
On the other hand,  by \eqref{connection1},
\beq \Lambda_{\omega_g}\smo D^{\prime\prime}\left(D^{\prime,h_0}H \cdot H^{-1}\right)
=\Lambda_{\omega_g}\smo D^{\prime\prime}\left(\p^{h_0}H \cdot H^{-1}-(\theta_{h_0}^{\star}\cdot H-H\cdot \theta_{h_0}^{\star}) H^{-1}\right). \eeq
For degree reasons,  one has
\beq \Lambda_{\omega_g}\smo \bp\left(\p^{h_0}H \cdot H^{-1}-\theta_{h_0}^{\star}+H\cdot \theta_{h_0}^{\star}\cdot H^{-1}\right)=\Lambda_{\omega_g}\smo\bp\left( (\p^{h_0} H) \cdot H^{-1} \right), \eeq
and by Lemma \ref{keycomputation},
\begin{align*}
&\Lambda_{\omega_g}\smo \theta\left(\p^{h_0}H \cdot H^{-1}-\theta_{h_0}^{\star}+H\cdot \theta_{h_0}^{\star}\cdot H^{-1}\right)\\
=&-\Lambda_{\omega_g}\smo\left( \theta\cdot\left(\p^{h_0}H \cdot H^{-1}-\theta_{h_0}^{\star}+H\cdot \theta_{h_0}^{\star}\cdot H^{-1}\right)\right)\\
&-\Lambda_{\omega_g}\smo\left(\left(\p^{h_0}H \cdot H^{-1}-\theta_{h_0}^{\star}+H\cdot \theta_{h_0}^{\star}\cdot H^{-1}\right)\cdot\theta\right)
\\=&\Lambda_{\omega_g}\smo\bigg(\theta \cdot  \theta_{h_0}^{\star}+ \theta_{h_0}^{\star}\cdot \theta-\theta \cdot  H\cdot \theta_{h_0}^{\star}\cdot H^{-1}- H\cdot \theta_{h_0}^{\star}\cdot H^{-1}\cdot \theta \bigg).
\end{align*}
Since $D''=\bp+\theta$,  we obtain
\begin{align*} & \Lambda_{\omega_g}\smo D^{\prime\prime}\left(D^{\prime,h_0}H \cdot H^{-1}\right)-\Lambda_{\omega_g}\smo\bp\left( (\p^{h_0} H) \cdot H^{-1} \right)\\=&\Lambda_{\omega_g}\smo\bigg(\theta \cdot  \theta_{h_0}^{\star}+ \theta_{h_0}^{\star}\cdot \theta-\theta \cdot  H\cdot \theta_{h_0}^{\star}\cdot H^{-1}- H\cdot \theta_{h_0}^{\star}\cdot H^{-1}\cdot \theta \bigg).
\end{align*}
By comparing with \eqref{difference4}, we derive \eqref{curvaturedifference}. \eproof

\vskip 1\baselineskip

\section{Comparison theorems for Higgs bundles}

In this section,  we establish Theorem \ref{main3},  which provides a comparison result for Higgs bundles. As applications, we obtain the uniqueness of solutions to the prescribed Hermitian-Yang-Mills tensor equation for Higgs bundles.
\btheorem   \label{CompareThmHiggs}
Let $(E,\theta)$ be a Higgs bundle over a compact Hermitian manifold $(M,\omega_g)$. Suppose that $h$ and $h_0$ are Hermitian metrics on $E$. If the Hermitian-Yang-Mills tensors of their Higgs connections satisfy $
\Lambda_{\omega_g}\left(\smo R^{D^{h_0}}\right) > 0 $ and \beq \Lambda_{\omega_g}\left(\smo R^{D^h}\right) \leq \Lambda_{\omega_g}\left(\smo R^{D^{h_0}}\right)  \eeq as Hermitian tensors in $ \Gamma(M,E^*\otimes \bar E{}^*) $, then $ h \leq h_0 $.
\etheorem

\noindent As an application of Theorem \ref{CompareThmHiggs}, one has

\bcorollary Let $(E,\theta)$ be a Higgs bundle over a compact Hermitian manifold $(M,\omega_g)$.  If  $h$ and $h_0$ are  Hermitian metrics on $E$ such that \beq \Lambda_{\omega_g}\left(\smo R^{D^h}\right)= \Lambda_{\omega_g}\left(\smo R^{D^{h_0}}\right) >0 \eeq as Hermitian tensors in $ \Gamma(M,E^*\otimes \bar E{}^*) $, then $ h = h_0 $.
\ecorollary

\noindent As another application of Theorem \ref{CompareThmHiggs}, one yields
\bcorollary\label{C2estimate} Let $(E,\theta)$ be a Higgs bundle over a compact Hermitian manifold $(M,\omega_g)$. Suppose that $h$ and $h_0$ are two Hermitian metrics on $E$. If there exists some $\lambda>0$ such that  $
\Lambda_{\omega_g}\left(\smo R^{D^{h_0}}\right) > 0 $ and \beq \Lambda_{\omega_g}\left(\smo R^{D^h} \right)\leq \lambda \cdot \Lambda_{\omega_g}\left(\smo R^{D^{h_0}}\right)  \eeq as Hermitian tensors in $ \Gamma(M,E^*\otimes \bar E{}^*) $, then $ h \leq \lambda h_0 $.
\ecorollary
\bproof A straightforward computation shows that the $(1,1)$-component of the curvature tensor $\Theta^{D^{h_0}}$ is
\beq \Theta^{1,1}:= \p^{h_0}\bp+\bp\p^{h_0}+\theta^\star_{h_0}\circ \theta+\theta\circ\theta_{h_0}^\star \in \Gamma(M,\Lambda^{1,1}T^*M\ts E^*\ts E).\eeq
This $(1,1)$-component is clearly invariant under scaling transformations of $h_0$.
Let $R^{D^{\lambda h_0}}$ be the curvature tensor of the Higgs bundle $(E, \lambda h_0, \theta)$. Then
\beq  \Lambda_{\omega_g}\left(\smo R^{D^{\lambda h_0}}\right)=\lambda \cdot \Lambda_{\omega_g}\left(\smo R^{D^{h_0}}\right) \in \Gamma(M, E^*\ts \bar E^*). \eeq
Hence, the result follows by Theorem \ref{CompareThmHiggs}.
\eproof

\noindent The proof of Theorem \ref{CompareThmHiggs} adopts a strategy analogous to that developed in the proof of \cite[Theorem~2]{WYY26+}. The primary differences and difficulties lie in two key aspects: the need to estimate the torsion of the Hermitian manifold, and the presence of the additional nonzero Higgs field $\theta$.
\bproof[Proof of Theorem \ref{CompareThmHiggs}] Let $\Omega =K^{h_0}= \Lambda_{\omega_g}\left(\smo \Theta^{D^{h_0}}\right)  \in \Gamma(M,E^*\otimes E)$ be the Hermitian-Yang-Mills tensor of $(E,h_0,\theta)$. For simplicity, we define $\Phi\in \Gamma(M,E^*\ts E)$ by the relation
\beq \Phi :=S^h\cdot h_0^{-1}=K^h\cdot H=  \left(\Lambda_{\omega_g}\left(\smo \Theta^{D^{h}}\right)\right)\cdot H,\eeq
where $H=h\cdot h_0^{-1} \in \Gamma(M,E^*\ts E)$.
By \eqref{curvaturedifference}, we obtain the identity
\beq \label{Currelation}\Phi\cdot H^{-1} = \Omega + \smo[\Lambda_{\omega_g},D''](D'^{,h_0}H \cdot H^{-1}). \eeq
This equation admits the following equivalent form:
\beq \label{mainidentity2} \Omega \cdot (H^{-1} - \mathrm{Id}_E) = \smo[\Lambda_{\omega_g},D''](D'^{,h_0}H \cdot H^{-1}) + (\Omega - \Phi) \cdot H^{-1}. \eeq
We claim that \beq H \leq \operatorname{Id}_E \label{keycompare}\eeq  with respect to $h_0$ which  implies $h \leq h_0$.
To prove \eqref{keycompare},  we define $\kappa: M \rightarrow \mathbb{R}$ as
\beq
\kappa(x):=\sup _{v \in E_x, v \neq 0} \frac{ h_0(Hv,v) }{h_0(v,v)} .
\eeq
The maximal eigenvalue of $H$ with respect to $h_0$ is
\beq
\Lambda:=\sup _{x \in M} \kappa(x) .
\eeq
We argue by contradiction and assume that \eqref{keycompare} does not hold. In this case, one has \beq  \Lambda > 1. \eeq

\noindent
For any $ m \geq 1 $,  pairing with $ H^m $ in \eqref{mainidentity2}  yields the following result:
\beq \label{PairingHm} \left(\Omega \cdot (H^{-1}-\mathrm{Id}_E),H^m \right)_{h_0} =\left(\smo[\Lambda_{\omega_g},D''](D'^{,h_0}H \cdot H^{-1}), H^m\right)_{h_0}+ \left((\Omega - \Phi) \cdot H^{-1}, H^m \right)_{h_0}. \eeq
The curvature condition asserts that $ \Omega - \Phi \geq 0 $ with respect to $h_0$.  By Lemma \ref{linearalgebra},
\beq \left((\Omega - \Phi) \cdot H^{-1}, H^m \right)_{h_0} = \left(\Omega-\Phi, H^{m-1}\right)_{h_0} \geq 0. \label{positivity} \eeq
On the other hand, by Bochner-Kodaira formula \eqref{bkformula} for Higgs bundles, one has
\be \left(\smo[\Lambda_{\omega_g},D''](D'^{,h_0}H \cdot H^{-1}), H^m\right)_{h_0} &=& \left(\left((D'^{,h_0})^*+\tau^*\right)(D'^{,h_0}H \cdot H^{-1}), H^m\right)_{h_0} \\
& = & \left(D'^{,h_0}H \cdot H^{-1}, D'^{,h_0}H^m+\tau  H^m \right)_{g, h_0}. \ee
By Lemma \ref{Lebniz},  $D'^{,h_0}H^m=\sum\limits_{p=0}^{m-1} H^p\cdot D'^{,h_0} H \cdot H^{m-1-p} $. By Lemma \ref{linearalgebra} again, one concludes
\beq\left\langle D'^{,h_0}H \cdot H^{-1}, D'^{,h_0}H^m \right\rangle _{g, h_0}=\left \langle D'^{,h_0}H, \sum\limits_{p=0}^{m-1} H^p\cdot D'^{,h_0} H \cdot H^{m-2-p}  \right\rangle_{g, h_0}.\eeq
Moreover, they are real functions.
We claim that, for any $ 0 \leq p \leq m-1 $, the following pointwise estimate holds:
\beq  \left\langle D'^{,h_0}H , H^p \cdot D'^{,h_0}H \cdot H^{m-2-p} \right\rangle_{g, h_0} + \frac{|\tau|_g^2}{4m^2}\left\langle \mathrm{Id}_E, H^m\right\rangle_{h_0} \geq -\frac{1}{m} \mathrm{Re}\left\langle D'^{,h_0}H, \tau  H^{m-1} \right\rangle_{g, h_0}. \label{pmatrixinequality} \eeq
Indeed, since $H$ is positive definite, we set
\beq A_p := H^{p/2} \cdot D'^{,h_0} H\cdot H^{(m-2-p)/2} + \frac{1}{2m} \tau  H^{m/2}.  \eeq
By Lemma \ref{linearalgebra}, one can see clearly that the estimate \eqref{pmatrixinequality} is equivalent to the fact that \beq \langle A_p,A_p\rangle_{h_0} \geq 0.\eeq  Summing \eqref{pmatrixinequality} over $ p=0$ to $p=m-1$ yields
\beq \sum_{p=0}^{m-1} \left\langle D'^{,h_0}H , H^p \cdot D'^{,h_0}H \cdot H^{m-2-p} \right\rangle_{g,h_0} + \frac{|\tau|_g^2}{4m}\left\langle \mathrm{Id}_E, H^m \right\rangle_{h_0} \geq -\mathrm{Re}\left\langle D'^{,h_0}H, \tau  H^{m-1} \right\rangle_{g, h_0}, \eeq
 and so
\beq \left\langle D'^{,h_0}H \cdot H^{-1}, D'^{,h_0}H^m\right\rangle_{g, h_0} + \mathrm{Re}\left\langle D'^{,h_0}H \cdot H^{-1}, \tau  H^m \right\rangle_{g, h_0} \geq -\frac{|\tau|_g^2}{4m}\langle\mathrm{Id}_E, H^m \rangle_{h_0}. \eeq
In particular,
\beq \label{torsionestimate}\left(\smo[\Lambda_{\omega_g},D''](D'^{,h_0}H\cdot H^{-1}), H^m\right)_{h_0} \geq -\frac{1}{4m}\|\tau\|^2_{L^\infty(M,\omega_g)} (\mathrm{Id}_E,H^m)_{h_0}. \eeq
Let $ K $ and $ k_0 $ be the maximal and the minimal eigenvalues of $ \Omega $ with respect to $h_0$ respectively, i.e.,
\beq K = \sup_{x\in M} \sup_{v \in E_x, v\neq 0} \frac{\nm{\Omega v}_{h_0}}{\nm{v}_{h_0}} > 0, \quad k_0 = \inf_{x\in M} \inf_{v \in E_x, v\neq 0} \frac{\nm{\Omega v}_{h_0}}{\nm{v}_{h_0}} > 0. \eeq
By \eqref{PairingHm}, \eqref{positivity} and \eqref{torsionestimate}, one has
\beq \label{Postiveofpairing}\left(\Omega, H^{m-1}-H^m \right)_{h_0} + \frac{C}{m}(\mathrm{Id}_E,H^m)_{h_0} = \left(\Omega \cdot (H^{-1}-\mathrm{Id}_E),H^m \right)_{h_0}+ \frac{C}{m}(\mathrm{Id}_E,H^m)_{h_0} \geq 0, \eeq
where $ C =\frac{1}{4}\|\tau\|^2_{L^\infty(M,\omega_g)} > 0 $ is a constant. \\

We shall employ a different method to bound the left-hand side of \eqref{Postiveofpairing} and reach a contradiction.  The following facts hold.
\bd
\item For any point $ x \in M $, one has
\beq  \left\langle \Omega,H^{m-1}-H^m \right\rangle_{h_0} \leq rK, \label{III}\eeq
where $ r = \mathrm{rank}(E) $ is rank of $ E $.
\item For any $ x \in M $ satisfying $ \kappa(x) \geq 1 $, one has
\beq \left\langle \Omega, H^{m-1}-H^m \right\rangle_{h_0} \leq rK + k_0\left(\kappa(x)^{m-1} - \kappa(x)^m\right). \label{IV}\eeq
\ed
Indeed, let $ f $ be the function that
\beq f(x) = x^{m-1}-x^m,\quad x\in(0,\Lambda]. \eeq
Let $ 0 < \lambda_1 \leq \cdots \leq \lambda_r \leq \Lambda $ be the eigenvalues of $ H $ at point $ x \in M $. The eigenvalues of $ f(H) $ are $ \{f(\lambda_k)\}_{k=1,2\cdots,r} $. It is clear that
$ f(x) \leq 1$ for $x \in (0,\Lambda]$.
In particular,
$ f(H) \leq \mathrm{Id}_E$.
Therefore,
\beq \left\langle \Omega, H^{m-1}-H^m \right\rangle_{h_0} \leq\left\langle \Omega, \mathrm{Id}_E \right\rangle_{h_0} \leq rK. \eeq
This establishes $(1)$.  For $(2)$,  if  $ x \in M $ satisfies $ \kappa(x) \geq 1 $, then one has
\beq f(\lambda_r) = f(\kappa(x)) = \kappa(x)^{m-1} - \kappa(x)^m. \eeq
Suppose that $ v_x\in E_x$ is  an eigenvector with $h_0(v_x,v_x)=1$  corresponding to $ \lambda_r = \kappa(x) $, then
\beq f(H) \leq \mathrm{Id}_E + \left(\kappa(x)^{m-1} - \kappa(x)^m\right) \; v_x^* \otimes v_x . \eeq
In particular, one has
\be  \left\langle \Omega, H^{m-1}-H^m \right\rangle_{h_0}(x) & \leq & \left\langle \Omega, \mathrm{Id}_E \right\rangle_{h_0} + \left(\kappa(x)^{m-1} - \kappa(x)^m\right) \left\langle \Omega, v_x^* \otimes v_x \right\rangle_{h_0} \\
& \leq & rK + k_0\left(\kappa(x)^{m-1} - \kappa(x)^m\right).  \ee
This gives \eqref{IV}. By  claims $(1)$ and $(2)$, one obtains
\beq \label{Pairingestimate1} \left( \Omega,H^{m-1}-H^m\right)_{h_0} \leq rK\mathrm{Vol}(M,\omega_g) + k_0\int_{\kappa(x)\geq 1} \left( \kappa(x)^{m-1} - \kappa(x)^m\right) \cdot \frac{\omega_g^n}{n!}. \eeq
and
\beq \label{Pairingestimate2}\frac{C}{m}(\mathrm{Id}_E,H^m)_{h_0}=\frac{C}{m}\int_M \tr H^m \frac{\omega_g^n}{n!} \leq \frac{Cr}{m}\mathrm{Vol}(M,\omega_g) + \frac{Cr}{m}\int_{\kappa(x)\geq 1} \kappa(x)^m \cdot \frac{\omega_g}{n!}. \eeq
By \eqref{Postiveofpairing}, \eqref{Pairingestimate1} and \eqref{Pairingestimate2}, we obtain
\beq \label{Postiveestimate} \left(\frac{Cr}{m}+rK\right)\mathrm{Vol}(M,\omega_g) + \int_{\kappa(x)\geq 1} \left( k_0\kappa(x)^{m-1} - \left(k_0-\frac{Cr}{m}\right)\kappa(x)^m\right) \cdot \frac{\omega_g^n}{n!} \geq 0.\eeq
Fix an arbitrary $\lambda \in (1,\Lambda)$, there exist $ m_0 \in \mathbb{N} $ and $ \delta> 0 $ such that for any $ m \geq m_0 $,
\beq \lambda\left(k_0-\frac{Cr}{m}\right) > k_0+\delta. \eeq
We set $ \lambda_1 = (\lambda+\Lambda)/2 $. \bd \item[(I)]  If $ \kappa(x) \geq  \lambda_1>\lambda $ and $ m \geq m_0 $, one has
\beq  k_0\kappa(x)^{m-1} - \left(k_0-\frac{Cr}{m}\right)\kappa(x)^m \leq \kappa(x)^{m-1}\left( k_0 - \lambda\left(k_0-\frac{Cr}{m}\right)\right) < -\delta\lambda_1^{m-1}. \eeq
\item[(II)] If $ \lambda \leq  \kappa(x) < \lambda_1 $ and $ m \geq m_0 $, one has
\beq  k_0\kappa(x)^{m-1} - \left(k_0-\frac{Cr}{m}\right)\kappa(x)^m \leq \kappa(x)^{m-1}\left( k_0 - \lambda\left(k_0-\frac{Cr}{m}\right)\right) < 0. \eeq
\item[(III)] If $ 1 \leq  \kappa(x) < \lambda $ and $ m \geq m_0 $, one has
\beq k_0\kappa(x)^{m-1} - \left(k_0-\frac{Cr}{m}\right)\kappa(x)^m \leq k_0\lambda^{m-1}. \eeq
\ed
Therefore, we obtain
\beq \label{Intergralestimate} \int_{\kappa(x)\geq 1} \left( k_0\kappa(x)^{m-1} - \left(k_0-\frac{Cr}{m}\right)\kappa(x)^m\right) \cdot \frac{\omega_g^n}{n!} \leq -\delta\mu\lambda_1^{m-1} + k_0\lambda^{m-1}\mathrm{Vol}(M,\omega_g), \eeq
where $ \mu > 0 $ is the measure of the set $\{x\in M\ | \kappa(x)\geq  \lambda_1\} $. By \eqref{Postiveestimate} and \eqref{Intergralestimate}, one has
\beq \left(\frac{Cr}{m}+rK + k_0\lambda^{m-1}\right)\mathrm{Vol}(M,\omega_g) - \delta\mu\lambda_1^{m-1} \geq 0, \eeq
and so
\beq \label{Contradiction}\lambda_1^{-(m-1)}\left(\frac{Cr}{m}+rK + k_0\lambda^{m-1}\right)\mathrm{Vol}(M,\omega_g) - \delta\mu \geq 0. \eeq
Note that the constants $\lambda_1, \lambda, C, r, K, k_0, \delta$ and $\mu$ are independent of $m$.
Moreover, since $ \lambda_1 > \lambda > 1 $,  when $ m \rightarrow +\infty $, by \eqref{Contradiction} we obtain $ -\delta\mu \geq 0 $. This is a contradiction. The proof is completed.
\eproof

\vskip 1\baselineskip

\section{A prior estimates}

Let $(M,\omega_g)$ be a compact Hermitian manifold and $(E,\theta)$ be a Higgs bundle.  Suppose that there exists a  Hermitian metric $h_0$ on $E$ such that
\beq \Om:=\Lambda_{\omega_g}\left(\sq \Theta^{D^{h_0}}\right)>0\eeq
with respect to $h_0$.
If $h$ is another  Hermitian metric,  we set \beq \Phi:=  \left(\Lambda_{\omega_g}\left(\smo \Theta^{D^{h}}\right)\right)\cdot H,\eeq
where $H=h\cdot h_0^{-1}$.
The formula \eqref{curvaturedifference} can be  written into
\beq
\Phi - \Omega \cdot H = \Lambda_{\om_g}\smo D''(D'^{,h_0}H \cdot H^{-1}) \cdot H. \label{curvaturedifference1}
\eeq
For degree reasons, it is clear that
\beq
\Lambda_{\om_g}\smo D''(D'^{,h_0}H \cdot H^{-1})=\Lambda_{\om_g}\smo\bp(\p^{h_0}H \cdot H^{-1})+\Lambda_{\om_g}\smo\theta(\theta_{h_0}^\star(H)\cdot H^{-1}).
\eeq
For simplicity, we set  \beq F_{\theta}(H):=\Lambda_{\om_g}\smo\theta(\theta_{h_0}^\star(H)\cdot H^{-1}). \eeq Then the equation \eqref{curvaturedifference1}  takes the form
\beq
\Phi-\Omega\cdot H-F_{\theta}(H)\cdot H=\Lambda_{\om_g}\smo\bp(\p^{h_0}H \cdot H^{-1}) \cdot H.\label{higgs equation}
\eeq

\noindent
The main result of this section is following uniform $C^1$-estimate for $H$.

\btheorem \label{ThmC1estimate}  Let $(E,\theta)$ be a Higgs bundle over a compact Hermitian manifold $(M,\omega_g)$. Suppose that $h_0$ and $h$ are two Hermitian metrics on $E$, and $\Lambda_{\omega_g}\left(\sq R^{D^{h_0}}\right)>0$.
If there exists a constant  $C_1\geq 1 $ such that
\beq C_1^{-1}h_0 \leq \Lambda_{\omega_g}\left(\sq R^{D^h}\right)\leq C_1 h_0 \eeq
as Hermitian tensors in $\Gamma(M,E^*\ts \bar E^*)$,
then  the following  $C^1$-estimate holds for $H=h\cdot h_0^{-1}$:
\beq  |H|_{C^1(M,\omega_g, h_0)}\leq C_{2},\eeq
where $C_{2}$ depends on $M,\omega_g, C_1, h_0,\theta$ and an upper bound of $|\Phi|_{C^1(M,\omega_g, h_0)}$.
\etheorem

\bremark The proof of Theorem \ref{ThmC1estimate} follows a variant of the strategy employed in \cite[Theorem~4.1]{WYY26+}, with careful attention paid to estimating both the torsion of the Hermitian manifold and the Higgs fields.
\eremark
\noindent We begin with the $C^0$-estimate which follows by Corollary  \ref{C2estimate}.
\bproposition
\label{C0estimate} Let $(E,\theta)$ be a Higgs bundle over a compact Hermitian manifold $(M,\omega_g)$. Suppose that $h_0$ and $h$ are two Hermitian metrics on $E$, and $\Lambda_{\omega_g}\left(\sq R^{D^{h_0}}\right)>0$.
If there exists a constant  $C_1\geq 1 $ such that
\beq C_1^{-1} h_0 \leq \Lambda_{\omega_g}\left(\smo R^{D^{h}}\right) \leq C_1 h_0,\eeq
then there exists a constant $ C_3 = C_3(M,\omega_g,h_0,\theta,C_1)  $ such that
\beq C_3^{-1} h_0 \leq h \leq C_3h_0.\eeq
\eproposition

\noindent The following estimate is established in \cite[Proposition~4.3]{WYY26+}
\bproposition
\label{W12estimate}   Let $E$ be a holomorphic vector bundle over  a compact Hermitian manifold $(M,\omega_g)$. Suppose that $h_0$ and $h$ are two Hermitian metrics on $E$, and $H=h\cdot h_0^{-1}$.
Suppose that there exists a constant $C_4\geq 1$ such that
\beq C_4^{-1}h_0\leq h\leq C_4h_0. \eeq
Then for any $A\in\Gamma(M,E^*\ts E)$,
\beq C^{-1}_4 h_0(A\cdot H^{-1}, A)\leq h(A\cdot H^{-1}, A\cdot H^{-1}) \leq C_4 h_0(A\cdot H^{-1}, A)\label{estimate1}\eeq
and
\beq C_4^{-2} h_0(A,A) \leq h(A\cdot H^{-1}, A\cdot H^{-1}) \leq C_4^2 h_0(A,A).\label{c15}\eeq
In particular,
\beq C_4^{-2} h_0(\p^{h_0} H, \p^{h_0} H) \leq h(\p^{h_0} H\cdot H^{-1}, \p^{h_0} H\cdot H^{-1}) \leq C_4^2 h_0(\p^{h_0} H, \p^{h_0} H).\eeq
Moreover,
\beq \mathrm{tr}_E\left(\Lambda_{\omega_g}\left(\sq \p^{h_0} H\cdot H^{-1}\cdot \bp H\right)\right)\geq C_4^{-1} h\left(\p^{h_0}H\cdot H^{-1}, \p^{h_0}H\cdot H^{-1}\right). \label{estimate2}\eeq
\eproposition

\bproposition \label{trEH} Let $(E,\theta)$ be a Higgs bundle over a compact Hermitian manifold $(M,\omega_g)$. Suppose that $h_0$ and $h$ are two Hermitian metrics on $E$, and $\Lambda_{\omega_g}\left(\sq R^{D^{h_0}}\right)>0$.
 If  there exists a constant $C_4\geq 1$ such that
\beq C_4^{-1}h_0\leq h\leq C_4h_0,\label{C0assume} \eeq
then
\beq \Delta_{\C}\mathrm{tr}_E H  \geq  C_4^{-1}|\p^{h_0}H \cdot H^{-1}|^2_h -\mathrm{tr}_E\Phi-C_5, \eeq
where $C_5=C_5(M,\om_g,h_0,\theta,C_4)$ and $\Delta_\C =\mathrm{tr}_{\omega_g}\sq \p\bp$ on functions.
\eproposition

\bproof A straightforward computation gives
\begin{eqnarray}
\Lambda_{\omega_g}\smo\bar\p(\p^{h_0}H \cdot H^{-1}) \cdot H  \nonumber& = & \Lambda_{\omega_g}\smo\bar\p(\p^{h_0}H \cdot H^{-1}) \cdot H \\
& = &  \Lambda_{\omega_g}\smo\bar\p\p^{h_0} H + \Lambda_{\omega_g}\smo(\p^{h_0}H \cdot H^{-1} \cdot \bar\p H). \label{calabic31}
\end{eqnarray}
Moreover, $\mathrm{tr}_E H=h_0(H, \mathrm{Id}_E)$ and so
\be \Delta_\C \mathrm{tr}_E H=\Lambda_{\omega_g} \sq \p\bp h_0(H, \mathrm{Id}_E)=-\sq \Lambda_{\omega_g} h_0(\bp \p^{h_0} H, \mathrm{Id}_E)=-\tr_E(\Lambda_{\omega_g}\smo\bar\p\p^{h_0} H). \ee
On the other hand, by \eqref{higgs equation}, one has
\beq \Phi-\Omega\cdot H-F_{\theta}(H)\cdot H=\Lambda_{\om_g}\smo \bp(\p^{h_0}H \cdot H^{-1}) \cdot H. \eeq
By taking the trace  of the equation \eqref{calabic31}, one obtains
\be \Delta_{\C}\mathrm{tr}_E H & =& -\tr_E(\Lambda_{\omega_g}\smo\bar\p\p^{h_0} H) \\ &=& \mathrm{tr}_E\left(\Lambda_{\omega_g}\smo(\p^{h_0}H \cdot H^{-1} \cdot \bar\p H)\right) - \mathrm{tr}_E(\Lambda_{\omega_g}\smo\bp (\p^{h_0}H \cdot H^{-1}) \cdot H ) \\
& \geq & C_4^{-1}|\p^{h_0}H \cdot H^{-1}|^2_h + \mathrm{tr}_E(\Omega \cdot H+F_\theta(H)\cdot H - \Phi) \\
& \geq &C_4^{-1}|\p^{h_0}H \cdot H^{-1}|^2_h+\tr_E(F_\theta(H)\cdot H) -\mathrm{tr}_E\Phi , \ee
where the first inequality follows from \eqref{estimate2} and the second inequality holds since $ \Omega >0$ with respect to $h_0$. It remains to estimate $\tr_E(F_\theta(H)\cdot H)$.  By Lemma \ref{keycomputation},
\begin{align*}
F_\theta(H)\cdot H=&\Lambda_{\om_g}\smo\theta(\theta_{h_0}^\star (H)\cdot H^{-1})\cdot H\\
=&\Lambda_{\om_g}\smo\theta(H\cdot \theta_{h_0}^\star \cdot H^{-1}-\theta_{h_0}^\star )\cdot H\\
=&\Lambda_{\om_g}\smo\left( -H\cdot \theta_{h_0}^\star \cdot H^{-1}\cdot\theta\cdot H-\theta\cdot H\cdot\theta_{h_0}^\star +\theta_{h_0}^\star \cdot \theta\cdot H+\theta\cdot\theta_{h_0}^\star \cdot H  \right).
\end{align*}
From \eqref{C0assume} one deduces that
\beq
C_4^{-1}\mathrm{Id}_E\leq H\leq C_4\mathrm{Id}_E.
\eeq
Hence
\beq
\mathrm{tr}_E(F_\theta(H)\cdot H)\geq- C_5(M,\om_g,h_0,\theta,C_4).
\eeq
This completes the proof.
\eproof

Let $ \sE: = T^{*1,0}M \otimes E^* \otimes E $ be the Hermitian holomorphic vector bundle with the Hermitian metric induced by $ g $ on $ T^{*1,0}M $ and $ h $ on $ E $.  Let  $ T $ be the tensor
\beq T := \p^{h_0}H \cdot H^{-1}  \in \Omega^{1,0}(M,E^*\otimes E) \simeq \Gamma(M,\sE). \eeq
We derive the following estimate:
\bproposition
\label{DeltaS} Let $(E,\theta)$ be a Higgs bundle over a compact Hermitian manifold $(M,\omega_g)$. Suppose that $h_0$ and $h$ are two Hermitian metrics on $E$.   If  there exists a constant $C_4\geq 1$ such that
\beq C_4^{-1}h_0\leq h\leq C_4h_0, \eeq
then
one has
\beq \Delta_{\C} |T|_h^2\geq - C_{21}\left(|T|^2_h+1\right) \eeq
where $C_{21}$ depends on $M,\omega_g, C_4, h_0,\theta$ and an upper bound of $|\Phi|_{C^1(M,\omega_g, h_0)}$.
\eproposition

\bproof For simplicity, we set
\beq
\Phi'=\left(\Lambda_{\om_g}\left(\smo \Theta^{h}\right)\right)\cdot H,
\eeq
where $\Theta^h$ is the Chern curvature tensor of $(E,h)$. By formulas \eqref{relation2} and \eqref{chernhiggs},
\begin{eqnarray}
\Phi-\Phi'&=&\left(\Lambda_{\om_g}\left(\smo \Theta^{D^h}\right)-\Lambda_{\om_g}\left(\smo \Theta^h\right)  \right)\cdot H\nonumber\\
&=&-\Lambda_{\om_g}\smo(\theta\cdot \theta_{h}^\star +\theta_{h}^\star \cdot \theta)\cdot {H}\nonumber\\
&=&-\Lambda_{\om_g}\smo(\theta\cdot H\cdot \theta_{h_0}^\star +H\cdot \theta_{h_0}^\star \cdot H^{-1}\cdot\theta\cdot H). \label{difference5}
\end{eqnarray}
Since $H$ is uniformly bounded, one obtains immediately that
\beq
|\Phi-\Phi'|_{h}\leq C_6(M,\omega_g, C_4, h_0, \theta).
\eeq
Moreover, since
\beq
T=\p^{h_0}H\cdot H^{-1}=-H\cdot \p^{h_0}H^{-1},
\eeq
One has
\beq C_4^{-1}|T|_h\leq |\p^{h_0}H|_h\leq C_4|T|_h, \quad C_4^{-1}|T|_h\leq |\p^{h_0}H^{-1}|_h\leq C_4|T|_h.
\eeq
By \eqref{difference5},  one concludes that
\beq
|\p^{h_0}\Phi-\p^{h_0}\Phi'|_h\leq C_7(|T|_h+1),
\eeq
where $C_7=C_7(M,\omega_g, C_4, h_0, \theta)$.
In particular, the following estimates hold:
\beq
|\Phi'|_h\leq C_8,\quad |\p^{h_0}\Phi'|_h\leq C_{9}(1+|T|_h).\label{chernbound}
\eeq
Here $C_8$, $C_{9}$ depend on $M,\omega_g, C_4, h_0, \theta$ and an upper bound of  $|\Phi|_{C^1(M,\omega_g, h_0)}$.

\vone
The Bochner-Kodaira formula establishes that
\beq \Delta_{\C} |T|_h^2 = |\p_\sE T|_h^2 + |\bar\p T|_h^2 + 2\mathrm{Re}\left(\left\langle \mathrm{tr}_{\omega_g}\left(\smo \p_\sE\bar\p_\sE T\right), T \right\rangle_\sE \right) - \mathrm{Ric}^\sE(T,T),\label{c30} \eeq
where $\mathrm{Ric}^\sE=\Lambda_{\omega_g} \left(\sq \Theta^{\sE}\right)$ is the Hermitian-Yang-Mills tensor of $\sE$.
It is clear that
\beq\mathrm{tr}_{\omega_g}\left(\smo \p_\sE\bar\p_\sE T\right) = \left(g^{i\bar j} \nabla^h_i \nabla^h_{\bar j} T_{k\alpha}^\beta\right) \; dz^k \otimes e^\alpha \otimes e_\beta,\label{c31}\eeq
where $\nabla^h$ is the Chern connection on $\sE$ induced by $g$ and $h$. Recall that   \beq T= T_{k\alpha}^\beta \; dz^k \otimes e^\alpha \otimes e_\beta=\left( (\Gamma^h)_{i\alpha}^\beta - (\Gamma^{h_0})_{i\alpha}^\beta \right) dz^i \otimes e^\alpha \otimes e_\beta.\eeq
Hence,
\begin{eqnarray} -g^{i\bar j} \nabla^h_i \nabla^h_{\bar j} T_{k\alpha}^\beta \nonumber& = & g^{i\bar j}\nabla^h_i \left((R^h)_{k\bar j\alpha}^\beta - (R^{h_0})_{k\bar j\alpha}^\beta \right) \\
& = & g^{i\bar j}\nabla^h_i (R^h)_{k\bar j\alpha}^\beta - g^{i\bar j}(\nabla_i^h - \nabla_i^{h_0})(R^{h_0})_{k\bar j\alpha}^\beta - g^{i\bar j}\nabla^{h_0}_i (R^{h_0})_{k\bar j\alpha}^\beta, \label{c32}
\end{eqnarray}
where $\nabla^{h_0}$ is the Chern connection induced by $h_0$ and $\omega_g$.
The torsion tensor of the Hermitian metric $g$ is denoted by
\beq
\Theta_{ij}^k = g^{k\bar\ell}\left(\frac{\p g_{j\bar\ell}}{\p z^i}-\frac{\p g_{i\bar\ell}}{\p z^j}\right).
\eeq
By using the torsion Bianchi identity on Hermitian manifolds, one has
\beq \nabla^{h}_i (R^{h})_{k\bar j\alpha}^\beta=\nabla^{h}_k (R^{h})_{i\bar j\alpha}^\beta+\Theta_{ki}^m(R^{h})_{m\bar j\alpha}^\beta, \eeq
and so we get
\be -g^{i\bar j} \nabla^h_i \nabla^h_{\bar j} T_{k\alpha}^\beta
&=&g^{i\bar j}\nabla^h_k (R^h)_{i\bar j\alpha}^\beta+\Theta_{ki}^m(R^{h})_{m\bar j\alpha}^\beta - g^{i\bar j}(\nabla_i^h - \nabla_i^{h_0})(R^{h_0})_{k\bar j\alpha}^\beta - g^{i\bar j}\nabla^{h_0}_i (R^{h_0})_{k\bar j\alpha}^\beta\\
& = & \nabla^h_k (\Phi'\cdot H^{-1})_\alpha^\beta +\Theta_{ki}^m(R^{h})_{m\bar j\alpha}^\beta- g^{i\bar j} (\nabla_i^h - \nabla_i^{h_0})(R^{h_0})_{k\bar j\alpha}^\beta -g^{i\bar j}\nabla^{h_0}_i (R^{h_0})_{k\bar j\alpha}^\beta.
\ee
Moreover,
\beq
\nabla^h_k (\Phi' \cdot H^{-1})_\alpha^\beta = \left(\nabla^h_k \left(\Phi'\right)_\alpha^\gamma\right) (H^{-1})_\gamma^\beta + \left(\Phi'\right)_\alpha^\gamma \nabla^h_k(H^{-1})_\gamma^\beta. \eeq
Since $h_0$ and $h$ are equivalent,
\beq |H^{-1}|_h\leq \sqrt{r} C_4\eeq
and
\beq \nabla^h_k \left(\Phi'\right)_\alpha^\gamma = \nabla^{h_0}_k \left(\Phi'\right)_\alpha^\gamma + (\nabla^h_k-\nabla^{h_0}_k)\left(\Phi'\right)_\alpha^\gamma, \eeq
we conclude that
\beq |\left(\nabla^h_k \left(\Phi'\right)_\alpha^\gamma\right) (H^{-1})_\gamma^\beta|_h \leq C_{10}|\p^{h_0}\Phi'|_h+C_{11}|T|_h |\Phi'|_h.\eeq
Similarly,
\beq \nabla_k^h(H^{-1})^\alpha_\gamma = \nabla_k^h(h^{\alpha\bar\beta}h_{0,\gamma\bar\beta}) = h^{\alpha\bar\beta} \nabla_k^h h_{0,\gamma\beta} = h^{\alpha\bar\beta} \left(\nabla_k^h - \nabla_k^{h_0}\right) h_{0,\gamma\bar \beta}. \eeq
This implies
\beq |\nabla^h_k (\Phi' \cdot H^{-1})_\alpha^\beta|_h\leq C_{12}|\p^{h_0}\Phi'|_h+C_{13}|T|_h |\Phi'|_h \label{c33}\eeq
Since $\Theta$ depends only on $\om_g$,
\beq
\left|\Theta_{ki}^m(R^{h}-R^{h_0})_{m\bar j\alpha}^\beta\right|_h= \left|\Theta_{ki}^m(\bp T)_{\bar jm\alpha}^\beta\right|_h\leq C_{14}|\bp T|_h.
\eeq
In particular,
\beq
\left|\Theta_{ki}^m(R^{h})_{m\bar j\alpha}^\beta\right|_h\leq \left|\Theta_{ki}^m(R^{h}-R^{h_0})_{m\bar j\alpha}^\beta\right|_h+\left|\Theta_{ki}^m(R^{h_0})_{m\bar j\alpha}^\beta\right|_h\leq C_{14}|\bp T|_h+C_{15}.\label{torsion}
\eeq
It is easy to see that
\beq \left|(\nabla_i^h - \nabla_i^{h_0})(R^{h_0})_{k\bar j\alpha}^\beta\right|_h\leq C_{16} |T|_h \label{c34}\eeq
and
\beq \left|g^{i\bar j}\nabla^{h_0}_k (R^{h_0})_{i\bar j\alpha}^\beta\right|_h\leq C_{17}. \label{c35}\eeq
Here the constants $C_{10}$, $C_{11}$, $\cdots$,  $C_{17}$  depend on $M,\omega_g, h_0$ and $ C_4$.
By \eqref{c31}, \eqref{c32}, \eqref{c33}, \eqref{torsion}, \eqref{c34} and \eqref{c35},  we deduce that
\beq \left|\mathrm{tr}_{\omega_g}\left(\smo \p_\sE\bar\p_\sE T\right)\right|_h \leq C_{12}|\p^{h_0}\Phi'|_h+C_{13}|T|_h |\Phi'|_h+C_{14}|\bp T|_h+C_{16} |T|_h+C_{15}+C_{17}. \label{c36}\eeq
Note that $ \sE= T^{*1,0}M \otimes E^* \otimes E $ and the Hermitian-Yang-Mills tensor $\mathrm{Ric}^\sE$ is a linear combination of the Hermitian-Yang-Mills tensors  $\mathrm{Ric^{(2)}(\omega_g)}$ and $\Lambda_{\om_g}\left(\smo \Theta^{h}\right)$, i.e.,
\beq \mathrm{Ric}^\sE=- \mathrm{\mathrm{Ric^{(2)}(\omega_g)}} \ts \mathrm{Id}_{E^*\ts E}+\mathrm{Id}_{T^{*1,0}M}\ts \left(-\left(\Lambda_{\om_g}\left(\smo \Theta^{h}\right)\right)^t\ts \mathrm{Id}_E+\mathrm{Id}_{E^*}\ts \Lambda_{\om_g}\left(\smo \Theta^{h}\right)\right).\eeq
Since $h$ and $h_0$ are equivalent and $\Phi'=\Lambda_{\om_g}\left(\smo \Theta^{h}\right)\cdot H$,  the following estimate holds
\beq  \mathrm{Ric}^\sE(T,T)\leq C_{18} |T|^2_h \left(1 +|\Phi'|_h\right).\label{c37}\eeq
By using the Cauchy-Schwarz inequality,
\beq
|\bp T|_h^2-2C_{14}|\bp T|_h|T|_h\geq-C_{19}|T|_h^2\label{c38}
\eeq
By \eqref{c30}, \eqref{c36}, \eqref{c37} and \eqref{c38}, one has
\beq \Delta_{\C} |T|_h^2 \geq -C_{20}\left(|\p^{h_0}\Phi'|_h|T|_h+|T|^2_h |\Phi'|_h+|T|^2_h+|T|_h\right). \eeq
where $C_{18}$,  $C_{19}$ and $C_{20}$ depend on $M,\omega_g, C_4,  h_0$.
By \eqref{chernbound} and   the  Cauchy-Schwarz inequality,  we establish   \beq \Delta_{\C} |T|_h^2\geq - C_{21}\left(|T|^2_h+1\right) \eeq
where $ C_{21}$ depends on $M,\omega_g, C_4, \theta, h_0$ and an upper bound of $|\Phi|_{C^1(M,\omega_g, h_0)}$.
\eproof

\bproof[Proof of Theorem \ref{ThmC1estimate}] By Proposition \ref{C0estimate}, there exists $ C_4 = C_4(M,\omega_g,h_0,\theta,C_1) \geq  1 $ such that
\beq C_4^{-1} h_0 \leq h \leq C_4h_0.\eeq
Moreover,
\beq \mathrm{tr}_E\Phi\leq \sqrt{r} |\Phi|_{h_0} \eeq
By Proposition \ref{trEH} and Proposition \ref{DeltaS}, there exists $C_5=C_5(M,\om_g,h_0,\theta,C_4)$
\beq \Delta_{\C}\mathrm{tr}_E H  \geq  C_4^{-1}|T|^2_h -\mathrm{tr}_E\Phi-C_5\geq  C_4^{-1}|T|^2_h -\sqrt{r} |\Phi|_{h_0}-C_5, \eeq and
\beq \Delta_{g} |T|_h^2 \geq -C_{21}(|T|_h^2 + 1), \eeq

\noindent   For a large $ L > 0 $ such that $LC_4^{-1} - C_{21}>1$, one has
\beq \Delta_{\C} (|T|_h^2 + L\mathrm{tr}_E H) \geq (LC_4^{-1} - C_{21})|T|_h^2 - (LC_{22}+LC_5 + C_{21}), \eeq
where $C_{22}$ is an upper bound of $\sqrt{r} |\Phi|_{h_0}$.
At a maximal point $ p \in M $ of $ (|T|_h^2 + L\mathrm{tr}_E H) $,  we obtain the estimate
\beq |T|_h^2(p)  \leq C_{23}. \eeq
Since $\mathrm{tr}_EH \leq r C_4$, we get
\beq (|T|_h^2 + L\mathrm{tr}_E H)(p) \leq  C_{24}. \eeq
Therefore, for any point $ x \in M $, one has
\beq |T|^2_h(x) \leq (|T|_h^2 + L\mathrm{tr}_E H)(x)\leq (|T|_h^2 + L\mathrm{tr}_E H)(p) \leq C_{24}. \eeq
Since $T=\p^{h_0} H \cdot H^{-1}$, one has
\beq |\p^{h_0} H \cdot H^{-1}|^2_h \leq C_{24}. \eeq
By \eqref{c15}, one gets
\beq |\p^{h_0} H|_{h_0}^2 \leq C_4^2 C_{24}. \eeq
In summary, we establish the uniform $C^1$-estimate
$|H|_{C^1(M,\omega_g, h_0)}\leq C_{25}$ where $ C_{25}$  depends on $M,\omega_g,h_0, \theta, C_1$ and an upper bound of $ |\Phi|_{C^1(M,\omega_g, h_0)}$.
\eproof

\vskip 1\baselineskip

\section{Proofs of main theorems}
In this section we prove Theorem \ref{main1}:
\btheorem
\label{Higgsequation}
Let $(E,\theta)$ be a Higgs bundle over a compact Hermitian manifold $(M,\omega_g) $. Suppose that there exists a smooth Hermitian metric $ h_0 $ on $E$ such that $ \Lambda_{\omega_g}\left(\sq  R^{D^{h_0}}\right)>0$. Then for any Hermitian positive definite tensor $ P\in \Gamma(M,E^*\otimes \bar E^*) $, there exists a unique smooth Hermitian metric $ h $ on $E$ such that  \beq \Lambda_{\omega_g}\left(\sq R^{D^h}\right)=P.\eeq
\etheorem
\noindent We begin by proving the following result, which ensures the closeness property required in Theorem \ref{Higgsequation}.

\btheorem\label{Ricciclose} Let $(E,\theta)$ be a Higgs bundle over a compact Hermitian manifold $(M,\omega_g) $.  Suppose that $ \{h_m\}$ is a sequence of smooth Hermitian metrics on $E$ and   $ S^{h_m} $ are the corresponding Hermitian-Yang-Mills-Higgs tensors in $\Gamma(M, E^*\ts \bar E^*)$.     If
\beq \lim_{m}\|S^{h_m}-P\|_{C^\infty(M,\omega_g, h_0)}=0 \eeq
for some Hermitian positive definite $P\in  \Gamma(M,E^* \otimes \bar E^*) $, then there exists a unique smooth Hermitian metric $ h $ on $E$ such that
\beq \Lambda_{\omega_g}\left( \sq R^{D^h}\right)=P. \eeq
\etheorem
\bproof We follow the ideas in the proof of \cite[Theorem~5.3]{WYY26+}. We fix notations in the space $\Gamma(M,E^*\otimes E)$:
\beq H_m := h_m\cdot  h^{-1}_{0}, \quad \Phi_m := K^{h_m}\cdot H_m=S^{h_m}\cdot h_0^{-1}, \quad  P^{h_0}: = P\cdot h_0^{-1}.\eeq Hence,  $ P^{h_0} > 0 $ with respect to $h_0$.   Since $\lim_{m}\|S^{h_m}-P\|_{C^\infty(M, \omega_g,h_0)}=0 $, one has \beq \lim_m\|\Phi_m-\Phi\|_{C^\infty(M, \omega_g,h_0)}=0. \eeq In particular,  there exist uniform constants $ c_1, c_2 > 0 $ and $m_0>0$, depending on $M,\omega_g, h_0$ and $P$, such that for any $ m\geq m_0 $, the following estimates holds
\beq c_1^{-1}\mathrm{Id}_E \leq \Phi_m \leq c_1\mathrm{Id}_E ,\quad  \nm{\Phi_m}_{C^1(M,\omega_g, h_0)} \leq c_2. \label{bound1} \eeq
By Proposition \ref{C0estimate}, there exists $c_3=c_3(M,\omega_g, h_0, c_2)$ such that
\beq c^{-1}_3 \mathrm{Id}_E\leq H_m\leq c_3\mathrm{Id}_E.\label{close-1}\eeq
By Theorem \ref{ThmC1estimate}, there exists $ C_0 = C_0(M,\omega_g,h_0,r, c_1,c_2,\theta) > 0 $ such that,
\beq
\|H_m\|_{C^1(M,\om_g,h_0)}\leq C_0.\label{bound2}
\eeq
In the proof of Proposition \ref{trEH}, we establish the relation
\beq
\Lambda_{\omega_g}\smo\bar\p\p^{h_0} H_m=\Phi_m-\Omega\cdot H_m-F_{\theta}(H_m)\cdot H_m-\Lambda_{\omega_g}\smo(\p^{h_0}H_m \cdot H_m^{-1} \cdot \bar\p H_m),
\eeq
where $\Omega = \Lambda_{\omega_g}\left(\smo \Theta^{D^{h_0}}\right)  \in \Gamma(M,E^*\otimes E)$ and
\be &&F_{\theta}(H_m)=\Lambda_{\om_g}\smo\theta(\theta_{h_0}^\star(H_m)\cdot H_m^{-1})\\&=&\Lambda_{\om_g}\smo\left( -H_m\cdot \theta_{h_0}^\star \cdot H_m^{-1}\cdot\theta\cdot H_m-\theta\cdot H_m\cdot\theta_{h_0}^\star +\theta_{h_0}^\star \cdot \theta\cdot H_m+\theta\cdot\theta_{h_0}^\star \cdot H_m  \right).
\ee
If we set
\beq
W_m:=\Phi_m-\Omega\cdot H_m-F_{\theta}(H_m)\cdot H_m-\Lambda_{\omega_g}\smo(\p^{h_0}H_m \cdot H_m^{-1} \cdot \bar\p H_m),\eeq then the following elliptic equation holds
\beq \Lambda_{\omega_g}\left(\smo\bar\p\p^{h_0} H_m\right)=W_m.
\eeq
By using \eqref{bound1}, \eqref{close-1} and \eqref{bound2}, the following estimate holds
\beq
\|W_m\|_{C^0(M,\om_g,h_0)}\leq C_1.
\eeq
By  $W^{k,p}$ and $C^{k,\alpha}$-estimates for elliptic systems, the conclusion follows from a similar argument as in the proof of \cite[Theorem~5.3]{WYY26+}.
\eproof

We consider the openness property. Let us recall some notations. The subspace of smooth Hermitian sections in $\Gamma(M,E^*\otimes \bar E^*)$ is denoted by $\mathrm{Herm}(E)$, and the subspace of smooth Hermitian metrics in $\Gamma(M,E^*\otimes \bar E^*)$ is denoted by $\mathrm{Herm}^+(E)$. By Corollary \ref{hermitian1} there is a natural map
$G: \mathrm{Herm}^+(E) \> \mathrm{Herm}(E) $ given by the Hermitian-Yang-Mills-Higgs tensor:
\beq G(h) = \Lambda_{\omega_g}\left(\sq R^{D^h}\right)\in \Gamma(M,E^*\otimes \bar E^*), \eeq
Fix a smooth Hermitian metric $ h_0 $ on $ E $, we define
\beq \mathrm{Herm}(E,h_0) := \left\{ S \in \Gamma(M,E^*\otimes E)|\ S \text{ is $h_0$-Hermitian}\right\}, \eeq
and
\beq \mathrm{Herm}^+(E,h_0) := \left\{ S \in \mathrm{Herm}(E,h_0)|\ S > 0  \right\}. \eeq
It is clear that the map $F:\mathrm{Herm}(E)\> \mathrm{Herm}(E,h_0)$ given by \beq  F(h) = h \cdot h_0^{-1} \eeq  is an isomorphism and its restriction $F:\mathrm{Herm}^+(E) \> \mathrm{Herm}^+(E,h_0)$ is also an isomorphism.
There is an induced map $\tilde G: \mathrm{Herm}^+(E,h_0) \> \mathrm{Herm}(E,h_0)$ \quad \beq   \tilde G = F \circ G \circ F^{-1}. \eeq
According to \eqref{curvaturedifference}, the map $ \tilde G: \mathrm{Herm}^+(E,h_0) \> \mathrm{Herm}(E,h_0) $ is given by
\beq \tilde G(H) = \Lambda_{\omega_g}\left(\smo D^{\prime\prime}\left(D^{\prime,h_0}H \cdot H^{-1}\right) \cdot H\right) + \Omega \cdot H, \eeq
where $ \Omega =  \Lambda_{\omega_g}\left(\sq \Theta^{D^{h_0}}\right) \in \mathrm{Herm}^+(E,h_0) $.   A straightforward computation shows that  $$D'^{, h_0}\left( \mathrm{Id}_E\right)=0,$$ and so the linearization of $ \tilde G $ at the point $ \mathrm{Id}_E $ is $\sL:\mathrm{Herm}(E,h_0)\rightarrow \mathrm{Herm}(E,h_0)$:
\beq \mathscr{L}(\Psi) = \Lambda_{\omega_g}\left(\smo D^{\prime\prime}D^{\prime,h_0}\Psi\right) + \Omega \cdot \Psi. \eeq
We can extend the definition of $\sL$ to $\Gamma(M,E^*\ts E)$.
\blemma
For any $\Psi\in \Gamma(M,E^*\otimes E)$,
\beq
\sL(\Psi^\star)=(\sL(\Psi))^\star\label{linearhermitianproperty}
\eeq
\elemma
\bproof
Since $\sL$ maps $\mathrm{Herm}(E,h_0)$ to $\mathrm{Herm}(E,h_0)$, one  has
\beq
\sL(\Psi+\Psi^\star)=\left(\sL(\Psi+\Psi^\star)\right)^\star, \quad
\sL\left(\smo(\Psi-\Psi^\star)\right)=\sL\left(\smo(\Psi-\Psi^\star)\right)^\star.
\eeq
Or equivalently,
\beq
\sL(\Psi)+\sL(\Psi^\star)=(\sL(\Psi))^\star+(\sL(\Psi^\star))^\star,\quad
\sL(\Psi)-\sL(\Psi^\star)=(\sL(\Psi^\star))^\star-(\sL(\Psi))^\star.
\eeq
Hence,  \eqref{linearhermitianproperty} follows.
\eproof

\btheorem\label{injectivity}  The map $\mathscr{L}:\Gamma(M,E^*\otimes E) \> \Gamma(M, E^*\otimes E) $ given by
\beq \mathscr{L}(\Psi) = \Lambda_{\omega_g}\left(\smo D^{\prime\prime}D^{\prime,h_0}\Psi\right) + \Omega \cdot \Psi \eeq
is injective if  $ \Omega  \in \mathrm{Herm}^+(E,h_0) $.
\etheorem
\bproof
Assume that $\sL(\Psi)=0$ for some $\Psi\in \Gamma(M,E^*\otimes E)$. By  \eqref{linearhermitianproperty},
\beq
\sL(\Psi)=\sL(\Psi^\star)=0
\eeq
Since the Chern connection is metric compatible,
\beq
\Delta_{g}|\Psi|_{h_0}^2=|\p^{h_0}\Psi|_{h_0}^2+|\bp\Psi|_{h_0}^2-h_0\left( \Lambda_{\om_g}\smo\bp\p^{h_0}\Psi,\Psi \right)+h_0 \left(\Psi,\Lambda_{\om_g}\smo\p^{h_0}\bp\Psi \right).\label{injectivemain}
\eeq
Note that
\begin{align*}
0=\sL(\Psi)=&\Lambda_{\omega_g}\left(\smo D^{\prime\prime}D^{\prime,h_0}\Psi\right) + \Omega \cdot \Psi\\
=&\Lambda_{\om_g}\smo\bp\p^{h_0}\Psi+\Lambda_{\om_g}\smo \theta(\theta_{h_0}^\star(\Psi))+\Omega\cdot\Psi.
\end{align*}
Hence
\beq
-\Lambda_{\om_g}\smo\bp\p^{h_0}\Psi=\Lambda_{\om_g}\smo \theta(\theta_{h_0}^\star(\Psi))+\Omega\cdot\Psi.
\eeq
Moreover,
\begin{align}
\begin{split}
-h_0\left( \Lambda_{\om_g}\smo\bp\p^{h_0}\Psi,\Psi \right)=&h_0\left( \Lambda_{\om_g}\smo \theta(\theta_{h_0}^\star(\Psi)),\Psi \right)+h_0\left( \Om\cdot\Psi,\Psi \right)\\
=&-\Lambda_{\om_g}\smo h_0\left(\theta^\star_{h_0}(\Psi),\theta^\star_{h_0}(\Psi)\right)+h_0\left( \Om\cdot\Psi,\Psi \right)\\
=&|\theta^\star_{h_0}(\Psi)|_{h_0}^2+h_0\left(\Om\cdot\Psi,\Psi\right)\geq 0.
\end{split}
\label{injectivemain1}
\end{align}
On the other hand,  by \eqref{adjoint},
\beq
(\p^{h_0}\bp\Psi)^\star=\bp((\bp\Psi)^\star)=\bp\p^{h_0}\Psi^\star,
\eeq
and so
\beq
h_0 \left(\Psi,\Lambda_{\om_g}\smo\p^{h_0}\bp\Psi \right)=h_0 \left((\Lambda_{\om_g}\smo\p^{h_0}\bp\Psi)^\star,\Psi^\star \right)
=-h_0\left(\Lambda_{\om_g}\smo\bp\p^{h_0}\Psi^\star,\Psi^\star\right).
\eeq
Note that $\sL(\Psi^\star)=0$, similar to \eqref{injectivemain1} one has
\beq
h_0 \left(\Psi,\Lambda_{\om_g}\smo\p^{h_0}\bp\Psi \right)=|\theta^\star_{h_0}(\Psi^\star)|_{h_0}^2+h_0(\Om\cdot\Psi^\star, \Psi^\star)\geq 0.\label{injectivemain2}
\eeq
By \eqref{injectivemain}, \eqref{injectivemain1} and \eqref{injectivemain2}
\beq
\Delta_{g}|\Psi|_{h_0}^2=|\p^{h_0}\Psi|_{h_0}^2+|\bp\Psi|_{h_0}^2+h_0(\Om\cdot\Psi,\Psi)+h_0(\Om\cdot\Psi^\star, \Psi^\star)+|\theta^\star_{h_0}(\Psi)|_{h_0}^2+|\theta^\star_{h_0}(\Psi^\star)|_{h_0}^2.
\eeq
Since $\Om$ is positive definite, by Lemma \ref{linearalgebra}, at a maximal point $p\in M$ of $|\Psi|_{h_0}^2$,
\beq
0\geq \Delta_{g}|\Psi|_{h_0}^2(p)\geq h_0(\Om\cdot \Psi,\Psi)(p)\geq 0,
\eeq
hence $\Psi(p)=0$, which implies that $\Psi\equiv 0$. \eproof

\bcorollary
\label{hermitian}
Suppose that  $ \Psi \in \Gamma(M,E^*\otimes E) $ satisfies
\beq \sL(\Psi)=\Lambda_{\omega_g}\left(\smo D^{\prime\prime}D^{\prime,h_0}\Psi\right) + \Omega \cdot \Psi = S. \eeq
If  $ S \in \mathrm{Herm}(E,h_0) $, then $ \Psi \in \mathrm{Herm}(E,h_0) $.
\ecorollary
\bproof
By  \eqref{linearhermitianproperty},
$
\sL(\Psi^\star)=S^\star=S=\sL(\Psi)$.
By Theorem \ref{injectivity},
$
\Psi=\Psi^\star$.
\eproof

\btheorem
\label{SecRicciopen}
Let $ P \in \mathrm{Herm}^+(E) $. If there exists a Hermitian metric $ h_1\in \mathrm{Herm}^+(E) $ such that
\beq P = G(h_1)=\Lambda_{\omega_g}\left( \sq R^{D^{h_1}}\right)  \in \Gamma(M,E^*\otimes \bar E^*), \eeq
then there exist open neighborhoods $ U $ of $ P \in  \mathrm{Herm}^+(E) $ and $ V $ of $ h_1 \in \mathrm{Herm}^+(E) $ such that $ G: U\> V $ is an isomorphism.

\etheorem

\bproof We fix the center point $h_1\in \mathrm{Herm}^+(E)$ and so $F(h_1)=\mathrm{Id}_E$.
Since $\tilde G=F\circ G\circ F^{-1}$,
\beq d\tilde G_{\mathrm{Id}_E}=(dF)_{G(h_1)}\circ (dG)_{h_1}\circ (dF^{-1})_{F(h_1)}\eeq
Moreover, since $ F $ is an isomorphism,  we only need to verify that  $$ \mathscr{L} = d\tilde G_{\mathrm{Id}_E}:\mathrm{Herm}(E,h_1) \> \mathrm{Herm}(E,h_1) $$ is an isomorphism. Let $\Om_1=\Lambda_{\omega_g}\left( \sq \Theta^{D^{h_1}}\right)  \in \Gamma(M,E^*\otimes  E)$.  Recall that
\begin{align*}
\sL(\Psi)=&\Lambda_{\omega_g}\smo D^{\prime\prime}D^{\prime,h_1}\Psi + \Omega_1 \cdot \Psi\\
=&\Lambda_{\omega_g}\smo \bp\p^{h_1}\Psi+ \Lambda_{\omega_g}\smo \theta(\theta_{h_1}^\star(\Psi)) +\Om_1\cdot\Psi\\
=&\p^{*,h_1}\p^{h_1}\Psi+\tau^*\p^{h_1}\Psi+\Lambda_{\omega_g}\smo \theta(\theta_{h_1}^\star(\Psi)) +\Om_1\cdot\Psi
\end{align*}
is an elliptic operator.  For any positive integer $m$, let
\beq W^{m,2}(M, E^*\ts E)\eeq
be the space of $W^{m,2}$ sections of $E^*\ts E$ with respect to the metric $\omega_g$ and $h_1$.
We claim that for all $ k \geq 1 $,
\beq \mathscr{L}: W^{k+2,2}(M, E^*\ts E) \> W^{k,2}(M, E^*\ts E)\label{Lk+2}\eeq
are isomorphisms. We only need to show the injectivity and surjectivity.
\bd
\item(Injectivity) If there exists  some $ \Psi \in W^{k+2,2}(M, E^*\ts E)$ such that $ \mathscr{L}(\Psi) = 0$, by elliptic regularity theory, $\Psi$ is  smooth. By Theorem \ref{injectivity}, $\Psi=0$.
\item(Surjectivity)
Consider the  pairing
\beq  \mathscr B: W^{1,2}(M, E^*\ts E)\times W^{1,2}(M, E^*\ts E)\> \mathbb{C} \eeq defined by
\beq \mathscr B (\Psi,T) =(\mathscr{L}(\Psi), T)_{h_1}\eeq in the weak sense.
It is easy to see that
\beq
\mathscr B (\Psi,T) =(\p^{h_1}\Psi,\p^{h_1}T)+(\p^{h_1}\Psi,\tau T)+(\Lambda_{\om_g}\smo\theta(\theta_{h_1}^\star(\Psi)) +\Om_1\cdot\Psi,T).
\eeq
By Fredholm alternative theorem (e.g. \cite[p.323]{Eva10}) and the fact that $\sL$ is injective, we deduce that for any $\Phi\in W^{1,2}(M,E^*\ts E)$, there exists some $\Psi\in W^{1,2}(M,E^*\ts E)$ such that $\sL(\Psi)=\Phi$ in the weak sense. By standard regularity theory for elliptic PDEs (e.g. \cite[p.~329]{Eva10}), one concludes that $ \Psi \in W^{k+2,2}(M, E^*\ts E)$.
\ed
Hence, for any  $\Phi\in  \mathrm{Herm}(E,h_1)$, there exists some $\Psi\in \Gamma(M,E^*\ts E)$ such that $\sL(\Psi)=\Phi$. By Corollary \ref{hermitian}, $\Psi\in \mathrm{Herm}(E,h_1) $. This implies the surjectivity of the linearization map  $
\sL:\mathrm{Herm}(E,h_1) \> \mathrm{Herm}(E,h_1)$. Hence, Theorem \ref{SecRicciopen} holds by the implicit function theorem.
\eproof

\btheorem
\label{SecondRicciCk}
Let $(E,\theta)$ be a Higgs bundle over a compact Hermitian manifold $(M,\omega_g) $. Suppose that there exists a smooth Hermitian metric $ h_0 $ on $E$ such that $ \Lambda_{\omega_g}\left(\sq  R^{D^{h_0}}\right)>0$. Then for any integer $k>1$, $\alpha\in (0,1)$ and  any  $ P \in C^{k,\alpha}{\mathrm{Herm}^+(E)} $, there exists a unique Hermitian metric $ h \in C^{k+2,\alpha}{\mathrm{Herm}^+(E)} $  such that its Hermitian-Yang-Mills-Higgs tensor $S^h$  coincides with  $P$.
\etheorem

\bproof
The uniqueness of metric $h$ is established by comparison theorem since $h$ is $C^3$. For the existence of such a metric,  we define
\beq \mathscr{R} = \{ P \in C^{k,\alpha}{\mathrm{Herm}^+(E)}  |\  S^h=P \text{ for some}\  h \in C^{k+2,\alpha}{\mathrm{Herm}^+(E)} \} . \eeq
By a similar argument as in the proof of Theorem \ref{SecRicciopen}, one deduces that  $ \mathscr{R} $ is open in $ C^{k,\alpha}{\mathrm{Herm}^+(E)}  $.
By using similar arguments as in the proof of Theorem \ref{Ricciclose},  one can show that:
if  $ \{h_m\}$ is a sequence of $C^{k+2,\alpha}$-Hermitian metrics on $E$ and
\beq \lim_{m}\|S^{h_m}-P\|_{C^{k,\alpha}(M,\omega_g, h_0)}=0 \eeq
for some  $P\in C^{k,\alpha}\mathrm{Herm}^+(E) $, then there exists a unique $h\in C^{k+2,\alpha}{\mathrm{Herm}^+(E)} $ such that
$S^h=P$.   Hence, $ \mathscr{R} $ is  closed in $ C^{k,\alpha} {\mathrm{Herm}^+(E)} $. Since $ C^{k,\alpha} {\mathrm{Herm}^+(E)} $ is connected,  one deduces that $  \mathscr{R}= C^{k,\alpha} {\mathrm{Herm}^+(E)}$.
\eproof

\bproof[Proof of Theorem \ref{Higgsequation}]
For any $ P \in \mathrm{Herm}^+(E) $, by Theorem \ref{SecondRicciCk}, there exists a unique Hermitian metric $ h_k \in C^{k+2,\alpha} {\mathrm{Herm}^+(E)} $ such that
\beq S^{h_k} = P. \eeq
Moreover, since $h_k$ is unique, $  h_k $  are the same for all $k>1$ and we denote it by $h$. Hence, $ h $ is a smooth Hermitian metric on $E$ and $ S^h = P $.
\eproof

    We demonstrate that Theorem \ref{main1} holds in a  more general setting.
Let $h$ be a smooth Hermitian metric on the Higgs bundle $ (E,\theta) $. The minimal eigenvalue function  $ \kappa_h: M \> \mathbb{R} $  of the Hermitian-Yang-Mills-Higgs tensor $K^h=\Lambda_{\omega_g}\left(\sq \Theta^{D^h}\right)\in \Gamma(M,E^*\ts E)$ is defined as
\beq \kappa_h(x) = \inf_{0 \neq v \in E_x} \frac{h(K^hv,v)}{h(v,v)}. \eeq

\btheorem\label{main6}
Let $ (E,\theta) $ be a Higgs bundle over a compact Hermitian manifold $(M,\omega_g) $. Suppose that $\omega_g$ is Gauduchon, i.e., $\p\bp\omega_g^{n-1}=0$, and  there exists a Hermitian metric $ h_0 $ on $E$ such that
\beq \int_M \kappa_{h_0}(x) \omega_g^n > 0. \label{integral}\eeq Then for any Hermitian positive definite tensor $ P\in \Gamma(M,E^*\otimes \bar E^*) $,  there exists  a unique smooth Hermitian metric $ h $ on $E$ such that  \beq \Lambda_{\omega_g} \left(\sq R^{D^h}\right)=P.\eeq
\etheorem

\bproof By using the condition \eqref{integral}, there exists  some $\lambda_0>0$ such that
\beq \int_M \kappa_{h_0}(x) \omega_g^n = \lambda_0\int_M \omega_g^n. \eeq
By the Hodge theory for Gauduchon metrics (e.g. \cite{Gau84}, \cite{Yang19a}),  there exists some function $ f \in C^{2,\alpha}(M,\mathbb{R}) $ such that
\beq \kappa_{h_0}+ \Delta_{\C} f = \lambda_0. \eeq
We define a new metric $\tilde h = e^{-f}h_0 $.  By \eqref{curvaturedifference}, one has
\beq \Lambda_{\omega_g} \left(\smo \Theta^{D^{\tilde h}}\right) -\Lambda_{\omega_g} \left(\smo \Theta^{D^{h_0}}\right) = \Lambda_{\omega_g}\smo D^{\prime\prime}\left(D^{\prime,h_0}H \cdot H^{-1}\right). \eeq
where $H=\tilde h\cdot h_0^{-1}= e^{-f} \cdot \mathrm{Id}_E$. A straightforward computation shows
\beq \Lambda_{\omega_g}\smo D^{\prime\prime}\left(D^{\prime,h_0}H \cdot H^{-1}\right)  = -\smo\Lambda_{\omega_g} \bar\p\p f \cdot \mathrm{Id}_E = \Delta_{\C}f \cdot \mathrm{Id}_E. \eeq
Moreover, one has
\beq \kappa_{\tilde h}(x) = \inf_{0 \neq v \in E_x} \frac{{\tilde h}(K^{\tilde h}v,v)}{{\tilde h}(v,v)} = \kappa_{h_0}(x) + \Delta_{\C}f(x)=\lambda_0. \eeq
In particular, the Hermitian-Yang-Mills-Higgs tensor $\Lambda_{\omega_g}\left(\sq R^{D^{\tilde h}}\right)$ is positive-definite. By a perturbation trick, one can find a smooth Hermitian metric $h_1$  in a small neighborhood of $\tilde h$ such that $\Lambda_{\omega_g}\left(\sq R^{ D^{h_1}}\right)$ is positive-definite.   Now Theorem \ref{main6} follows by  Theorem \ref{main1}.\eproof

\noindent The method in the proof of Theorem \ref{main6} can be extended to establish several new results that generalize the findings in \cite{LY12, Yang18,  Yang19a, Yang19b, Yang20, Yang24, XYY24+}. For instance,
\bcorollary Let $(M,\omega_g)$ be a compact K\"ahler manifold. Suppose that the integral of the minimal eigenvalue of the Ricci curvature $\mathrm{Ric}(\omega_g)$ is positive, then there exists a Hermitian metric $h$ on $M$, such that
\beq \Lambda_{\omega_g}\left(\sq R^{h}\right)=g.\eeq
Moreover, one has $H_{\bp}^{p,0}(M)=0$ for $p=1,\cdots, \dim M$. In particular, $M$ is projective and rationally connected.
\ecorollary
\noindent
We  propose:
\begin{problem} Let $(M,\omega_g)$ be a compact K\"ahler manifold. \bd \item  Suppose that the integral of the minimal eigenvalue of the holomorphic sectional curvature is positive. Does this imply that $M$ is projective and rationally connected?
    \item  Suppose that the integral of the maximal eigenvalue of the holomorphic sectional curvature is negative (resp. nonpositive). Does this imply that $K_M$ is ample (resp. nef)? \ed
\end{problem}

\vskip 1\baselineskip

\section{Quantitative Chern number inequalities for Higgs bundles}

In this section, we prove Theorem \ref{main5}:
\btheorem Let $(M,\omega_g)$ be a compact K\"ahler manifold and $(E,\theta)$ be an integrable Higgs vector bundle of rank $r$ over $M$.  Suppose that there exists a smooth Hermitian metric $h_0$ on $E$ such that the Hermitian-Yang-Mills tensor of the Higgs connection satisfies
\beq  a\cdot h_0 \leq   \Lambda_{\omega_g}\left(\smo  R^{D^{h_0}}\right) \leq b \cdot h_0, \eeq
for some constants $a, b\in \R$. Then the following Chern number inequality holds
\beq \int_M \left( (r-1)c^2_1(E)-2rc_2(E) \right) \wedge \omega_g^{n-2} \leq \frac{r(r-1)\left(b - a\right)^2}{8\pi^2n^2}\int_M \omega_g^n.\label{CN} \eeq
\etheorem
\bproof
Since $(E,\theta)$ is integrable,
\beq
(D'')^2=\bp\circ \bp +\bp\circ \theta+\theta\circ \bp+\theta\circ \theta.
\eeq Since $\bp\theta=0$ and $\theta\wedge \theta=0$, $(D'')^2=0$.
One can also show $(D'^{,h_0})^2=0$.
Therefore for any Hermitian metric $h$ on $E$, \beq \Theta^{D^{h}}=D'^{,h}D''+D''D'^{,h}=\p^{h}\bp+\bp\p^{h}+\theta^\star_h\circ \theta+\theta\circ\theta_h^\star+\p^{h}\theta+\bp\theta_h^\star. \eeq  In particular,
\beq
\left(\Theta^{D^{h}}\right)^{(2,0)}=\p^{h}\theta,\quad
\left(\Theta^{D^{h}}\right)^{(0,2)}=\bp\theta^\star_{h}, \quad \left(\Theta^{D^{h}}\right)^{(1,1)}= \p^{h}\bp+\bp\p^{h}+\theta^\star_h\circ \theta+\theta\circ\theta_h^\star.
\eeq
For the sake of simplicity,  we denote  $R=\Theta^{D^h} $. By \eqref{adjoint}, we know
\beq R^{(2,0)}=\left(R^{(0,2)}\right)^\star. \label{curvatureadjoint}\eeq
By the Chern-Weil theory,  we have the following Chern forms representations:
\beq
c_1(E)=\frac{\smo}{2\pi}\tr_ER,\quad 2c_2(E)=-\frac{1}{4\pi^2}\left(\tr_E R\wedge \tr_ER-\tr_E\left( R\wedge R \right)\right).
\eeq
We consider the $(1,1)$-component of $R$ and define
\beq
\tilde c_1(E):=\frac{\smo}{2\pi}\tr_E R^{(1,1)},\quad 2\tilde c_2(E):=-\frac{1}{4\pi^2}\left(\tr_E R^{(1,1)}\wedge \tr_ER^{(1,1)}-\tr_E\left(R^{(1,1)}\wedge R^{(1,1)}\right)\right).
\eeq
Moreover, one has the following identities for the $(2,2)$-form components
\beq
\left(c_1^2(E)-\tilde c_1^2(E)\right)^{(2,2)}=-\frac{1}{2\pi^2}\tr_E\left(R^{(2,0)} \right)\wedge \tr_E\left(R^{(0,2)} \right), \eeq
and \beq
\left(c_2(E)-\tilde c_2(E)\right)^{(2,2)}=-\frac{1}{4\pi^2}\left(\tr_E\left(R^{(2,0)} \right)\wedge \tr_E\left(R^{(0,2)} \right)-\tr_E\left(R^{(2,0)}\wedge R^{(0,2)} \right)\right).
\eeq
In particular,
\begin{align*}
&\int_M\left( (r-1)c^2_1(E)-2rc_2(E) \right)\wedge \omega_g^{n-2}-\int_M\left( (r-1)\tilde c^2_1(E)-2r\tilde c_2(E) \right)\wedge \omega_g^{n-2}\\
=&\frac{1}{2\pi^2}\int_M\tr_E\left(R^{(2,0)} \right)\wedge \tr_E\left(R^{(0,2)} \right)\wedge\omega_g^{n-2}-\frac{r}{2\pi^2}\int_M\tr_E\left(R^{(2,0)}\wedge R^{(0,2)} \right)\wedge \omega_g^{n-2}\\
=&-\frac{r}{2\pi^2} \int_M \tr_E (\eta^{(2,0)}\wedge \eta^{(0,2)})\wedge \omega_g^{n-2},
\end{align*}
where
\beq \eta^{(2,0)}=R^{(2,0)}-\left(\tr_E R^{(2,0)}\right)\ts \frac{\mathrm{Id}_E}{r},\quad \eta^{(0,2)}=R^{(0,2)}-\left(\tr_E R^{(0,2)}\right)\ts \frac{\mathrm{Id}_E}{r}.\eeq
By \eqref{curvatureadjoint}, one can see clearly that
\beq \eta^{(2,0)}=\left(\eta^{(0,2)}\right)^\star.\eeq
In particular,
\beq  \int_M \tr_E (\eta^{(2,0)}\wedge \eta^{(0,2)})\wedge \omega_g^{n-2}\geq 0.\eeq
As a direct consequence,
\beq\int_M\left( (r-1)c^2_1(E)-2rc_2(E) \right)\wedge \omega_g^{n-2}-\int_M\left( (r-1)\tilde c^2_1(E)-2r\tilde c_2(E) \right)\wedge \omega_g^{n-2}\leq 0.
\eeq
We use the following notations to denote curvature tensors of the Higgs connection: $$R^{(1)}_{i\bar j} =h^{\alpha\bar\beta}R^{(1,1)}_{i\bar j \alpha\bar \beta}, \quad R^{(2)}_{\alpha\bar\beta}=g^{i\bar j}  R^{(1,1)}_{i\bar j \alpha\bar\beta } \qtq{and} s_h =g^{i\bar j} R^{(1)}_{i\bar j}.$$
The following expressions are well-known:
\beq \label{c1computation} \int_M \tilde  c^2_1(E) \wedge \omega_g^{n-2} = \frac{1}{4\pi^2n(n-1)}\int_M \left(s_h^2 - \left|\mathrm{Ric}^{(1)}\right|_g^2 \right) \omega_g^n, \eeq
and
\beq \label{c2computation} \int_M \tilde  c_2(E) \wedge \omega_g^{n-2} = \frac{1}{8\pi^2n(n-1)}\int_M \left(s_h^2 - \left|\mathrm{Ric}^{(1)}\right|_g^2 - \left| \mathrm{Ric}^{(2)}\right|^2_h + |R^{(1,1)}|_{g,h}^2\right) \omega_g^n. \eeq
We define $T\in \Gamma(M,\Lambda^{1,1}T^*M\ts E^*\ts \bar E^*)$:
\beq T_{i\bar j \alpha\bar\beta } = R^{(1,1)}_{i\bar j\alpha\bar\beta} - \frac{1}{n}g_{i\bar j}R^{(2)}_{\alpha\bar\beta} - \frac{1}{r}R^{(1)}_{i\bar j}h_{\alpha \bar\beta} + \frac{1}{nr}g_{i\bar j}h_{\alpha \bar\beta}s_h. \eeq
A straightforward computation shows that
\beq |T|^2=|R^{(1,1)}|^2 - \frac{1}{n}\left|\mathrm{Ric}^{(2)}\right|^2 - \frac{1}{r}\left|\mathrm{Ric}^{(1)}\right|^2 + \frac{1}{nr} s_h^2. \label{T2}\eeq
By using \eqref{T2}, one has
\be \int_M \left(2r\tilde c_2(E) - (r-1)\tilde c^2_1(E)\right) \wedge \omega_g^{n-2}= \frac{1}{4\pi^2n(n-1)}\int_M \left(r|T|^2 + \frac{n-1}{n}\left( s_h^2 - r\left|\mathrm{Ric}^{(2)}\right|^2\right) \right) \omega_g^n. \ee
    Let $\{\lambda_1,\cdots, \lambda_r\}$ be the eigenvalues of $ \Lambda_{\omega_g}\sq R^{D^h}=\sum R^{(2)}_{\alpha\bar\beta }e^\alpha\ts \bar e^\beta $  with respect to $ h $.  Since $ a h\leq \Lambda_{\omega_g}\sq R^{D^h} \leq bh $, it is obvious that $ a \leq \lambda_k \leq b$, and so
\beq s_h^2 - r\left|\mathrm{Ric}^{(2)}\right|^2 = \left(\sum_{k=1}^r \lambda_k\right)^2 - r\sum_{k=1}^r \lambda_k^2 =  -\sum_{i<j} (\lambda_i-\lambda_j)^2 \geq -\frac{r(r-1)}{2}(a-b)^2. \eeq
Finally,  we establish the Chern number inequality.
\eproof

    \vskip 1\baselineskip

    \section{The RC-positivity theory for Higgs bundles}
As noted in \cite{WYY26+}, the foundational motivation for this research program stems from the RC-positivity theory, developed over the series of studies  in \cite{Yang18,  Yang19a, Yang19b, Yang20, Yang24, XYY24+}. We begin by recalling the formal notation established for RC-positivity. Let $E$ be a holomorphic vector bundle over a
complex manifold $M$.  $E$  is called \emph{RC-positive}, if there exists a smooth Hermitian metric $h$ on $E$ such that for each $q\in M$ and each nonzero vector $v\in  E_q$, there exists {some} nonzero  vector $u\in T_qM$
such that \beq R^h(u,\bar u,v,\bar v)>0,\eeq
where $R^h$ is the curvature tensor of the Chern connection on $(E,h)$.

\noindent
We propose the notion of RC-positivity for Higgs bundles:
\bdefinition  Let $(E,\theta)$ be a Higgs bundle over a compact complex manifold $M$, i.e., $\bp\theta=0$.  It is called \emph{RC-positive}, if there exists a smooth Hermitian metric $h$ on $E$ such that for each $q\in M$ and each nonzero vector $v\in  E_q$, there exists {some} nonzero  vector $u\in T_qM$
such that \beq R^{D^h}(u,\bar u,v,\bar v)>0,\eeq
where $R^{D^h}$ is the  curvature tensor of the Higgs connection on $(E,h,\theta)$.

\edefinition

    \noindent We anticipate that the positivity constraint on the Hermitian-Yang-Mills-Higgs tensor $\Lambda_{\omega_g}\left(\sq R^{D^{h_0}}\right)$ stated in Theorem \ref{main1} is generalizable to the RC-positivity of $R^{D^{h_0}}$. For instance, Theorem \ref{main3} holds for RC-positive Higgs bundles.   A systematic study extending RC-positivity theory to the framework of Higgs bundles--incorporating both vanishing theorems and comparison theorems--is forthcoming.

    \vskip 1\baselineskip


\begin{thebibliography}{Kawamata22}



        \bibitem[Cal57]{Cal57}
        Calabi, E.  On K\"ahler manifolds with vanishing canonical class,  Algebraic geometry and topology Symposium in honor of S. Lefschetz, \emph{Princeton Univ. Press}, 1957, 78--89.


        \bibitem[Don85]{Don85}
        Donaldson, S. Anti self-dual Yang-Mills connections over complex algebraic surfaces and stable vector bundles,
         \emph{Proc. London Math. Soc.} \textbf{50} (1985), 1--26.

        \bibitem[Don87]{Don87}
        Donaldson, S.  Infinite determinants, stable bundles and curvature,
    \emph{Duke Math. J.} \textbf{54} (1987), 231--247.


        \bibitem[Eva10]{Eva10}
 Evans, L.-C.  {\it Partial differential equations}, second edition,  Amer. Math. Soc., Providence, RI, 2010.


        \bibitem[FY08]{FY08}  Fu, J.-X.  and Yau, S.-T. The theory of superstring with flux on non-K\"ahler manifolds and the complex Monge-Amp\`ere equation, \emph{J. Differential Geom.} {\textbf 78} (2008), no.~3, 369--428.

        \bibitem[Gau84]{Gau84}
     Gauduchon, P. La $1$-forme de torsion d'une vari\'et\'e{} hermitienne compacte, \emph{Math. Ann. }{\textbf 267} (1984), no.~4, 495--518.

        \bibitem[Hit87]{Hit87}
         Hitchin, N.-J. The self-duality equations on a Riemann surface,\emph{ Proc. London Math. Soc.} (3) {\textbf 55} (1987), no.~1, 59--126.







        \bibitem[LY12]{LY12}
        Liu, K.-F. and Yang, X.-K.  Geometry of Hermitian manifolds, \emph{Internat. J. Math.} {\textbf 23} (2012), no.~6, 1250055, 40 pp.





    \bibitem[NS65]{NS65}
    Narasimhan, M.-S.  and  Seshadri, C.-S.  Stable and unitary vector bundles on a compact Riemann surface,  \emph{Ann. of Math. }(2) {\textbf 82} (1965), 540--567.




        \bibitem[Sim88]{Sim88} Simpson, C.-T.  Constructing variations of Hodge structure using Yang-Mills theory and applications to uniformization, \emph{J. Amer. Math. Soc.} {\textbf 1} (1988), no.~4, 867--918.

    \bibitem[Sim92]{Sim92} Simpson, C.-T.  Higgs bundles and local systems,  \emph{Inst. Hautes \'Etudes Sci. Publ. Math.} No. \textbf{75} (1992), 5--95.

    \bibitem[Sim94a]{Sim94a} Simpson, C.-T. Moduli of representations of the fundamental group of a smooth projective variety. I, \emph{Inst. Hautes \'Etudes Sci. Publ. Math.} No. \textbf{79} (1994), 47--129.

    \bibitem[Sim94b]{Sim94b} Simpson, C.-T.
    Moduli of representations of the fundamental group of a smooth projective variety. II, Inst. Hautes \'Etudes Sci. Publ. Math. No. \textbf{80} (1994), 5--79.

        \bibitem[STW17]{STW17}
    Sz\'ekelyhidi, G. ; Tosatti, V. and Weinkove, B. Gauduchon metrics with prescribed volume form,  \emph{Acta Math.} \textbf{219}(2017), 181--211.



    \bibitem[TW10a]{TW10a} Tosatti, V.  and Weinkove, B.   Estimates for the complex Monge-Amp\`ere equation on Hermitian and balanced manifolds, \emph{Asian J. Math.} \textbf{14} (2010),  19--40.


    \bibitem[TW10b]{TW10b} Tosatti, V.  and Weinkove, B.  The complex Monge-Amp\`ere equation on compact Hermitian manifolds, \emph{J. Amer. Math. Soc.} \textbf{23} (2010), 1187--1195.

    \bibitem[WYY26+]{WYY26+} Wang, M.-W.; Yang, X.-K and Yau, S.-T.  Existence of Hermitian metrics  with  prescribed Hermitian-Yang-Mills tensors  I.  \emph{arXiv:2603.10611}

    \bibitem[XYY24+]{XYY24+} Xiong, Z.-Y.; Yang, X.-K and  Yau, S.-T.
    RC-positivity, Schwarz's lemma and comparison theorems. \emph{arXiv:2412.02553}

    \bibitem[Yang18]{Yang18}  Yang, X.-K.  RC-positivity, rational connectedness and Yau's conjecture, \emph{Camb. J. Math.} {\bf 6} (2018), no.~2, 183--212.

    \bibitem[Yang19a]{Yang19a}
    Yang, X.-K.  Scalar curvature on compact complex manifolds,  \emph{Trans. Amer. Math. Soc.} {\bf 371} (2019), no.~3, 2073--2087.

    \bibitem[Yang19b]{Yang19b}
    Yang, X.-K.  A partial converse to the Andreotti-Grauert theorem, \emph{Compos. Math.} {\bf 155} (2019), no.~1, 89--99;


    \bibitem[Yang20]{Yang20}
    Yang, X.-K.  RC-positive metrics on rationally connected manifolds, \emph{Forum Math. Sigma} {\bf 8} (2020), Paper No. e53, 19 pp.

    \bibitem[Yang24]{Yang24}
    Yang, X.-K.  Compact K\"ahler manifolds with quasi-positive second Chern-Ricci curvature, \emph{Comm. Anal. Geom.} {\bf 32} (2024), no.~10, 2717--2734.





        \bibitem[Yau78]{Yau78}
        Yau, S.-T.  On the Ricci curvature of a compact K\"ahler manifold and the complex Monge-Amp\`ere equation, \emph{Comm. Pure Appl. Math.} \textbf{31} (1978), 339--411.



        \bibitem[UY86]{UY86}
        Uhlenbeck, K. and Yau, S.-T.  On the existence of Hermitian-Yang-Mills connections in stable vector bundles, \emph{Comm. Pure Appl. Math.} \textbf{39} (1986), 257--293.


    \end{thebibliography}
\end{document}